%% file: DARP_MIP_TW.tex
\documentclass[11pt, oneside]{article}   	
\usepackage{geometry}                		
\usepackage{graphicx}				
\usepackage{amsmath}
\allowdisplaybreaks
\usepackage{amssymb}
\usepackage{xcolor}
\usepackage{longtable}
\usepackage{booktabs}
\usepackage{multirow}
\usepackage{colortbl}
\usepackage{natbib}
\usepackage{paralist}
\usepackage{pdflscape}
\usepackage{algpseudocode}
\usepackage{algorithm}
\usepackage{makecell}
\usepackage{tikz}
\usetikzlibrary{arrows,automata,shapes,calc,positioning, math, chains, intersections, arrows.meta}
\usetikzlibrary{decorations.markings, decorations.pathreplacing}

\usepackage{caption}
\usepackage{bbm}
\usepackage{amsthm}
\usepackage{mathtools}
\usepackage{cite}
\usepackage{subcaption}

\usepackage{setspace}
\usepackage{geometry}
\usepackage{hyperref}
\newcounter{Theo}
\newcounter{Algo}
\newcounter{Assu}

\theoremstyle{plain}

\title{A Branch-and-Cut algorithm for an urban dial-a-ride problem}

\usepackage{authblk}
\author[]{Arne Schulz\footnote{Corresponding author} \ and Christian Pfeiffer}
\affil[]{Universit\"at Hamburg, Institute of Operations Management, Moorweidenstra{\ss}e 18, 20148 Hamburg, Germany}
\affil[]{\textit{\{arne.schulz, christian.pfeiffer\}@uni-hamburg.de}}

\begin{document}
\maketitle
\begin{abstract}
The paper investigates a dial-a-ride problem focusing on the residents of large cities. These individuals have the opportunity to use a wide variety of transportation modes. Because of this, ridepooling providers have to solve the tradeoff between a high pooling rate and a small detour for customers to be competitive. We provide a Branch-and-Cut algorithm for this problem setting and introduce a new technique using information about already fixed paths to identify infeasible solutions ahead of time and to improve lower bounds on the arrival times at customer locations. By this, we are able to introduce additional valid inequalities to improve the search. We evaluate our procedure in an extensive computational study with up to 120 customers and ten vehicles. Our procedure finds significantly more optimal solutions and better lower and upper bounds in comparison with a mixed-integer programming formulation.
\end{abstract}

\section{Introduction}
The organization of passenger traffic is a major challenge for many cities. Citizens should be able to reach their destination quickly but congestion and pollution need to be avoided. While public transport with trains or buses generally reduces the amount of traffic greatly, it also offers little flexibility to the end user due to fixed time schedules and routes. On the other end of the spectrum, a private car or a taxi allows for more flexibility but requires a vehicle for each ride. Ridepooling services can solve this tradeoff by filling the gap between these two extremes and help to reduce traffic while providing highly flexible options to the customer. The resulting optimization problem of assigning customers to vehicles and determining the vehicles' routes is called dial-a-ride problem (DARP) in scientific literature. \\
\ \\
Ridepooling services focusing on the general public especially in large cities have already been implemented in practice. Examples are MOIA and ioki which both operate amongst others in Hamburg, Germany. \citet{TSSRSSR17} found that the potential for sharing rides and therefore for ridepooling services is massive in several large cities like New York, San Francisco, Singapore or Vienna while only generating detours of some minutes. \\
\ \\
To generate a sufficient market share, ridepooling providers need to be attractive enough to end users when compared with alternatives like a traditional taxi or different modes of transportation. Therefore, including the minimization of customer inconvenience in the objective function is of particular importance. However, due to time windows and a highly constrained vehicle capacity in practice, allowing customer requests to be rejected, i.e. not to be served, is necessary to obtain feasible solutions. Nevertheless, finding a solution where no customers are rejected is the primary objective. \\
\ \\
As ridepooling providers need to pool rides as efficiently as possible to operate economically while being sufficiently attractive for customers, strong exact solution approaches are necessary to solve the mentioned tradeoff optimally and evaluate ridepooling proposals in practice. Our contribution is twofold: We investigate the static DARP with ride and waiting time minimization first proposed by \citet{PS20} and extend the formulation by including time windows and allowing the rejection of customer requests. We solve this DARP variant with a Branch-and-Cut algorithm that exploits the problem structure to improve bounds and identify infeasible solutions ahead of time. We accomplish this by introducing additional valid inequalities and arc fixings both during the search and in the preprocessing. \\
\ \\
The remainder of the paper is structured as follows: In Section \ref{sec:literature}, we discuss the relevant literature before we define our problem setting formally and present a mixed-integer program including several valid inequalities in Section \ref{sec:problem}. In Section \ref{sec:BCA}, we describe the core of our Branch-and-Cut algorithm, the fixed path procedure. Afterwards, we present a computational study evaluating the Branch-and-Cut approach in Section \ref{sec:comp} and give a conclusion in Section \ref{sec:conclusion}.

\section{Literature review}\label{sec:literature}
There is a large amount of literature investigating dial-a-ride problems. Early works are by Psaraftis who developed exact algorithms to solve the problem (\citet{Psa80}, \citet{Psa83}). Other exact approaches involving Branch-and-Cut algorithms are presented by \citet{Cor06} and \citet{RCL07}. Here, the objective is to minimize the total routing costs while satisfying time windows and ride time and capacity constraints. Several model formulations and valid inequalities are proposed. In \citet{RC09}, the authors propose a Branch-and-Cut-and-Price algorithm with two variations for the pricing problem and some of the valid inequalities from the earlier papers. Their approach can outperform the previous Branch-and-Cut algorithms. Branch-and-Cut-and-Price is also used in \citet{GI15} who can improve on some of the earlier results. More recently, \citet{RF21} used partial tours (so called fragments) in a Branch-and-Cut framework which can decrease the solution times even further. Other works utilizing Branch-and-Cut approaches to solve the dial-a-ride or the similar pickup and delivery problem (PDP) include \citet{RR97}, \citet{LD04} and \citet{BCJ14}.\\
\ \\
\citet{BBM11} solved a PDP and considered two different objectives. The first is to minimize routing costs, the second is to first minimize the number of utilized vehicles and then the routing costs. The problem is solved via column generation. \citet{CMC10} utilized Branch-and-Cut with a Benders Decomposition to solve the PDP with transfers, i.e. in their problem a request can be picked up by one vehicle and then transferred to another vehicle that drives to the delivery location. Another problem variant that is solved with Branch-and-Cut-and-Price is the PDP with a heterogeneous vehicle fleet presented by \citet{QB15}. In this problem, the vehicles can be configured to change the capacity for the goods or passengers that are transported and different vehicle types offer different configuration options.\\
\ \\
The rejection of customers is also considered in \citet{RP06}. In their PDP, they allow for requests which are not served to be put into a request bank. Each of these requests is penalized in the objective function. The choice to not serve specific customers is also present in \citet{PPA15} where profits are considered. There it is advantageous to not transport a customer if the corresponding earnings do not exceed the costs of transport. The rejection of customers is similar to the selective DARP (\citet{RR18}, \citet{RF22}) in which the objective is to maximize the number of served requests which is our primary objective as well. \\
\ \\
Many of the papers that consider a minimization of customer inconvenience do so in multi-criteria objective functions where other objectives like routing costs are minimized as well. \citet{PDHG09} and \citet{MBCV17} minimized the total distance and the mean user ride time. \citet{DD04} examined a problem where the total distance, the excess ride time, and the total idle time are considered in the objective. In \citet{CLJD12}, the authors presented a DARP with three objectives: minimizing the number of vehicles, minimizing the duration of the vehicle tours, and minimizing the delays. \citet{PCLP13} investigated multiple factors that influence the customer satisfaction, one of which includes the squared ratio of actual to direct ride time. The measurement of the ride time relative to the direct route is also used by \citet{PS20}. Their objective function is almost the same as the one we use in this paper. The difference is that we do not require an additional factor to penalize the maximum detour, since we utilize time windows and additionally allow to reject customers as well. In contrast to many previous works, we do not minimize the total distance. As \citet{PS20} argued, ridepooling providers usually use electric cars which are cheaper to recharge such that driving costs are less important. On the other hand, generating a sufficient market share is crucial to reach an adequately high pooling rate, which is the reason why we focus on customer inconvenience. \\
\ \\
One paper that does not consider a dial-a-ride problem but shares similarities with our approach of successively building paths in a Branch-and-Bound algorithm is \citet{KHA85}. They considered a traveling salesman problem with pickup and delivery and used a Branch-and-Bound algorithm to solve the problem. The branching decision is whether or not to include an edge in the solution. However, they did not generate the lower bound of a node with the help of a linear programming relaxation but calculated their own lower bound based on reducing the distance matrix as shown in \citet{LMSK63}. \\
\ \\
For a more detailed discussion on literature of dial-a-ride problems, we refer to the surveys of \citet{PDH08} and \citet{HSKLPT18}.

\section{Problem definition}\label{sec:problem}
We consider a problem with $n$ customers each with an origin (pickup) and a destination (delivery) vertex. There are $K$ tours, each visiting a set of customers in a particular order. The tours start and end at a single depot. The vehicles are assumed to be homogeneous with a maximum passenger capacity $Q$. Each customer request is associated with a pickup location $i$ and a delivery location $n+i$ and has a number of passengers $q_i$. The visit of each pickup and delivery vertex has to happen within a hard time window given by $[e_i,l_i]$. We can interpret $e_i$ for a pickup location $i$ as the earliest point in time the customers want to start their ride and assume that our offer is not attractive enough for them if we are not able to serve them before $l_{n+i}$ whereat $n+i$ is the corresponding delivery location. Distances between locations $i$ and $j$ are represented by $t_{ij}$ and are assumed to be non-negative and to fulfill the triangle inequality. The complete notation is given in the following:
\begin{longtable*}{ll}
\label{tab:notation}\\
\mbox{\underline{sets:}}&\\
$P=\{1,\ldots, n\} $& \mbox{pickup locations} \\
$D=\{n+1,\ldots, 2n\} $&\mbox{delivery locations} \\
$i,j,h,k\in I$ & \mbox{locations}, $I= \{ 0, 2n+1 \}\cup P\cup D$, where $\{ 0, 2n+1 \}$ \\& represent the depot \\
& \\
\mbox{\underline{parameters:}}&\\
$t_{ij}$ & \mbox{travel time between locations $i$ and $j$} \\
$q_i$ & \mbox{demand at location $i$: $q_i>0$, $q_{i+n} = -q_i \quad\forall i\in P$}\\
$e_i$ & \mbox{start of time window for location $i\in P\cup D$}\\
$l_i$ & \mbox{end of time window for location $i\in P\cup D$}\\
$Q$ & \mbox{capacity of each vehicle}\\
$n$& \mbox{number of requests}\\
$K$&\mbox{number of tours}\\
$M_{ij}$&\mbox{sufficiently large constant}\\
$\Phi$ & \mbox{penalty factor for rejecting a person}\\
& \\
\mbox{\underline{variables:}}&\\
$X_{ij}$& \mbox{$= 1$ if a vehicle drives from location $i$ directly to location $j$,} \\ & 0 otherwise\\
$Q_i$& \mbox{number of passengers in the vehicle after leaving location $i$}\\
$B_i$ & \mbox{departure time from location $i$}\\
$Y_i$ & \mbox{$=1$ if the customer request of location $i$ is rejected, } 0 otherwise\\
\end{longtable*} \ \\
We minimize the customer inconvenience which we measure by the weighted relative detour proposed in \citet{PS20} and \citet{PS22}. The relative detour is calculated as the difference between the arrival at the delivery vertex and the sum of direct distance and earliest pickup time and put into relation to the direct distance. If customer $i$ is picked up at their earliest pickup time $e_i$ and driven straight to their delivery destination, their detour will be 0. This can be understood as a benchmark customers will compare their detour with, e.g. the use of a taxi or a private car. Due to the time window at the delivery vertex, each customer is effectively given an upper bound for their maximum detour. This is also the reason why we omit the term to penalize the maximum detour as presented in \citet{PS20}. \\
\ \\
We also allow the rejection of customers, i.e. to not serve the customer request. However, the rejection of a customer is penalized such that a solution with fewer rejected customers is always superior to a solution with more rejected customers. The penalty factor $\Phi$ is thus set to $\sum_{i\in P}q_i\frac{l_{i+n}-e_{i} - t_{i,i+n}}{t_{i,i+n}} + 1$ which is strictly larger than the worst possible assignment without rejection (if such an assignment exists). The ability to reject customers is beneficial in instances in which no feasible solution exists otherwise. In these instances, we can identify the customers who need to be rejected to find a solution where as many customers as possible are served in the best possible way. The penalty factor $\Phi$ is weighted with the number of passengers in a request $q_i$ such that the number of served people is maximized (which can differ from the number of served customer requests when there are requests with $q_i>1$). \\
\ \\
For every pickup vertex we assume w.l.o.g. that $e_i\geq t_{0i}$, i.e. every pickup vertex can be reached on time if visited first in a tour, and for every delivery vertex we assume that $e_{i+n}=e_i+t_{i,i+n}$ and $l_{i+n}\geq l_{i} + t_{i,i+n}$, i.e. a vehicle never has to wait at a delivery vertex. Consequently, if there is only a single customer in a tour, the contribution to the objective value for this customer is always 0.\\
\subsection{Mathematical model}\label{sec:model}
The following model formulation is based on the 2-index formulation by \citet{RCL07}. $\mathcal{S}$ is the set of all sets $S$ with $0 \in S$, $2n+1\notin S$ and at least one customer request $i$ where the delivery vertex is in $S$ but the pickup vertex is not, i.e., $\mathcal{S} = \{ S: 0 \in S \wedge 2n+1\notin S \wedge \exists i: (i\notin S \wedge n+i\in S) \}$. The model formulation follows as:

\begin{equation}
\min\quad \sum_{i\in P}q_i\frac{B_{i+n}-e_{i} - t_{i,i+n}}{t_{i,i+n}} + \sum_{i\in P}\Phi \cdot q_i\cdot Y_i\label{eq:objFunc}
\end{equation}
\begin{align}
&\sum_{i\in I}X_{ij} + Y_j =1\quad&&\forall j\in P\cup D\label{eq:degOut}\\ 
&\sum_{j\in I}X_{ij} + Y_i=1\quad&&\forall i\in P\cup D\label{eq:degIn}\\
&\sum_{j\in P}X_{0j}\leq K&&\label{eq:depot}\\
&\sum_{i,j\in S}X_{ij}\leq |S|-2\quad &&\forall S\in \mathcal{S}\label{eq:SEC}\\
&B_i+t_{ij}-M_{ij}(1-X_{ij}) \leq B_j &&\forall i\in P\cup D, j\in P\cup D\label{eq:time}\\
&Q_i+q_j-Q(1-X_{ij})\leq Q_j &&\forall i\in P\cup D, j\in P\cup D\label{eq:cap}\\
&Y_i = Y_{i+n}\quad &&\forall i\in P\label{eq:reject}\\
&e_i\leq B_i\leq l_i&&\forall i\in P\cup D \label{eq:B}\\
&\max\{0,q_i\}\leq Q_i\leq\min\{Q,Q+q_i\}&&\forall i\in P\cup D \label{eq:q}\\
&X_{ij}\in\{0,1\}&&\forall i\in I, j\in I\label{eq:xbin}\\
&Y_i\geq 0\quad&&\forall i\in P\cup D\label{eq:yNN}
\end{align}
The first part of the objective function \eqref{eq:objFunc} sums up the relative detours for every customer weighted with the number of passengers in the request. Since $B_{i+n}$ can never be smaller than $e_{i} + t_{i,i+n}$, the relative detour is always non-negative. The second part adds the penalty factor for every rejected passenger in the customer request. Equations \eqref{eq:degOut} and \eqref{eq:degIn} enforce that every pickup and delivery vertex is visited and left exactly once unless the customer request is rejected in which case the $Y_i$ variable will take the value 1. Restriction \eqref{eq:depot} makes sure that the depot is left at most $K$ times. Constraints \eqref{eq:SEC} are crucial in preventing violations of the precedence constraints. Without them, it would be possible for a customer's pickup and delivery vertices to be visited in the wrong order or even to be split across multiple tours. The Inequalities \eqref{eq:time} and \eqref{eq:cap} ensure that the time variables $B_i$ and capacity variables $Q_i$ are set correctly. $M_{ij}$ can be set to $\max(0,l_i+t_{ij}-e_j)$. Restrictions \eqref{eq:reject} enforce that the pickup and delivery vertex of a customer are either both rejected or not. In an actual implementation, the variables for the delivery vertices $Y_i$ with $i\in D$ can simply be omitted from the model and replaced with the corresponding pickup variable. The final Constraints \eqref{eq:B} -- \eqref{eq:yNN} restrict the variable domains. Due to the vertex degree Constraints \eqref{eq:degOut} and \eqref{eq:degIn}, variables $Y_i$ will always take binary values in an integer-feasible solution and do not need to be further constrained.\\
\ \\
Constraints \eqref{eq:SEC} are not added upfront but are checked for violations in every integer feasible solution like in \citet{RCL07}. This is accomplished by solving a maximum flow problem for every customer to see whether the pickup and delivery vertex are visited in the right order and in the same tour. Customers that are rejected in the current solution do not need to be checked. \\
\subsection{Preprocessing}\label{sec:preproc}
\begin{figure}[htbp]%
\begin{subfigure}[b]{0.6\textwidth}
\scalebox{0.8}{
\input{matching.pgf}
}
\subcaption{The three example sequences}
\label{fig:matchingI}
\end{subfigure}
\begin{subfigure}[b]{0.35\textwidth}
\scalebox{0.8}{
\input{matchingII.pgf}
}
\subcaption{The assignment problems}
\label{fig:matchingII}
\end{subfigure}
\caption{Example for the additional arc fixings. The left figure shows the three sequences of $i+n, j_1 \mbox{ or } j_2, \mbox{ and } k_1 \mbox{ and } k_2$ for which Conditions \eqref{eq:dominance} hold. In the right figure, the two resulting assignment problems are shown.}%
\label{fig:matching}%
\end{figure}
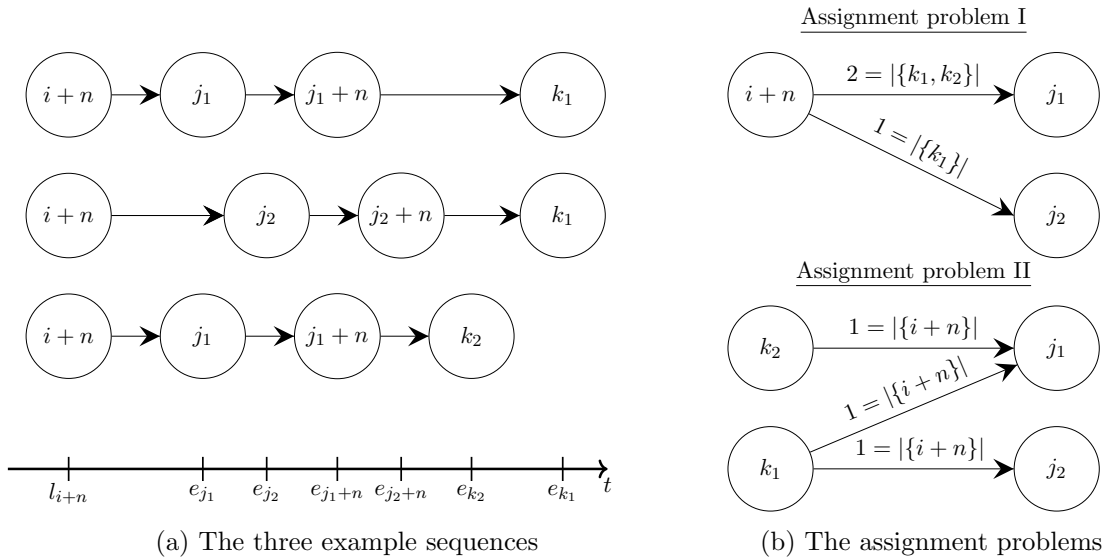

In the preprocessing, we identify incompatible customers like in \citet{Cor06} and fix arc variables to 0 where appropriate. Furthermore, we use additional inequalities to improve the bounds of the $Q_i$ and $B_i$ variables presented in \citet{RCL07} (Constraints (16) -- (19) in their paper).\\
\ \\
Moreover, we developed new preprocessing techniques presented in the following. We introduce new lower bounds for the $B_i$ variables. They are inspired by the dominance criteria for the dynamic program presented in \citet{PS20}. They function by identifying a pair of customers $i$ and $j$ where $i$ and $i+n$ can be visited before $j$ while not accumulating any kind of detour, i.e. we arrive at $j$ before the start of the time window $e_j$. Formally, if $e_j - t_{i,i+n}-t_{i+n,j}>e_i$, the following inequality is valid:
\begin{equation*}
\bar{B}_{ij} \leq B_i + (1-X_{ij})\cdot(\bar{B}_{ij} - e_i) \hspace{1cm} \forall i,j \in P
\end{equation*}
with $\bar{B}_{ij} = e_j-t_{i,i+n}-t_{i+n,j}$. Effectively, if we visit vertex $i$ at a time $B_i\leq \bar{B}_{ij}$, it is always better to visit the delivery vertex $i+n$ before heading to $j$. Additionally, when leaving the depot, this idea can be used to fix certain arc variables to 0. If we can leave the depot to visit customer $i$ and the delivery vertex $i+n$ and still arrive at customer $j$'s pickup location before or at the start of the time window $e_j$, it never makes sense to visit customer $j$ straight after the depot, so we set $X_{0j}$ to 0. If there are multiple customers $j$ for which this condition applies, we can only set one of the arcs to 0 because we can only visit the customer $i$ once. In such a case, we choose the customer with the lower pickup time window start $e_j$.\\
\ \\
Following the same idea, if there is a delivery vertex $i+n\in D$ and a pickup vertex $k\in P$ with a later time window, we can check if there is another customer $j\in P$ that we can visit perfectly in between the two vertices (i.e. the sequence $i+n,j,j+n,k$) which means that the vehicle arrives at or before the start of each time window (i.e. customer $j$ has an objective value contribution of 0). Additionally, the number of passengers of request $j$ cannot be larger than the maximum of the ones in requests $i$ and $k$: If $q_j\leq q_i$, the vehicle is guaranteed to have enough space to pick up request $j$ after leaving $i+n$. Should $q_i<q_j$ but $q_j\leq q_k$, this is not guaranteed, but in every case where request $j$ does not fit, the passengers of request $k$ would also not fit. This results in the following conditions:
\begin{align}
\begin{split}
&l_{i+n}+t_{i+n,j}\leq e_{j}\\ 
&e_{j+n}+t_{j+n,k}\leq e_k\\
&q_j\leq \max(q_i,q_k) 
\end{split}\label{eq:dominance}
\end{align}
If Conditions \eqref{eq:dominance} hold, arc $(i+n,k)$ can be set to 0, as visiting customer $j$ in between $i+n$ and $k$ is always superior to using the arc. While it is possible that an optimal solution uses the arc $(i+n,k)$, we could always replace it with the arc sequence $(i+n,j),(j,j+n),(j+n,k)$ and it would not make the solution worse, as both $j$ and $k$ are visited at their earliest possible time. Similar to the previous case each customer $j$ can only be used once, either for a vertex $i+n$ and all vertices $k$ for which the Conditions \eqref{eq:dominance} hold or for one vertex $k$ and all corresponding vertices $i+n$. Since we would like to fix as many arcs $(i+n,k)$ to 0 as possible, we are interested in finding an assignment of vertices $j$ such that this number is maximized. As a result, we solve two separate assignment problems (maximizing the total weights): First, we match vertices $i+n$ to vertices $j$ and the assignment ($i+n,j$) is weighted with the number of vertices $k$ for which Conditions \eqref{eq:dominance} apply. The second assignment problem does the reverse and matches vertices $k$ to $j$ and the matching ($k,j$) is weighted with the number of vertices $i+n$ for which Conditions \eqref{eq:dominance} hold. After both matching problems are solved, we choose the solution with the higher objective value (i.e. the higher number of fixed arcs $(i+n,k)$) and fix the corresponding arcs to 0. An example for this procedure is visualized in Figure \ref{fig:matching} where Figure \ref{fig:matchingI} shows the original sequences for which the conditions hold and Figure \ref{fig:matchingII} shows the two assignment problems. In this example, both problems have the same optimal objective value (assignment problem I: assign $j_1$ to $i+n$; assignment problem II: assign $j_1$ to $k_2$ and $j_2$ to $k_1$), but in general, this is not guaranteed. If we were to remove customer $j_2$ for example, the first assignment problem would still lead to an objective value of two, but the second assignment problem would only have an objective value of one. \\
\ \\
If we consider the solution of the first assignment problem in the figure, the assignment ($i+n,j_1$) is chosen in the optimal assignment. Thus, both arc $(i+n,k_1)$ as well as arc $(i+n,k_2)$ can be fixed to 0 since it is never worse to visit vertices $j_1$ and $j_1+n$ between them. Note that in general we can assign each vertex at most once because otherwise we might count an arc twice. In the first assignment problem for example, we can only assign either $j_1$ or $j_2$ to $i+n$ because otherwise there is an interaction which leads to a double counting of $k_1$ although we can fix the arc $(i+n,k_1)$ only once to zero.
\ \\
To further improve the bounding of the $B$ variables, we also add constraints of the form 
\begin{equation*}
B_i + t_{i,i+n}\cdot X_{i,i+n} + \sum_{j\in P\cup D\backslash\{i+n\}} (t_{ij}+t_{j,i+n})\cdot X_{ij}\leq B_{i+n}\quad \forall i\in P
\end{equation*}
which are valid since the sum of the $X$ variables can never exceed one. If there are request locations that are in the same spot, i.e. two locations $i,j\in P\cup D$ for which $t_{ij}=t_{ji}=0$, and if the time windows of these locations overlap, we can set the arc variable of the later to the earlier location to 0. If the earlier location is a pickup and the later a delivery, this last fixing is not used since it might violate the capacity constraint. Instead, if a pickup location $i$ and a delivery location $j$ have a distance $t_{ij}=t_{ji}=0$ and $i+n\neq j$, the variable $X_{ij}$ can be set to 0 (irrespective of the time windows). This is due to the assumption that we never have to wait at a delivery location and that the pickup location $j-n$ has to be visited prior. \\

\section{Branch-and-Cut algorithm}\label{sec:BCA}
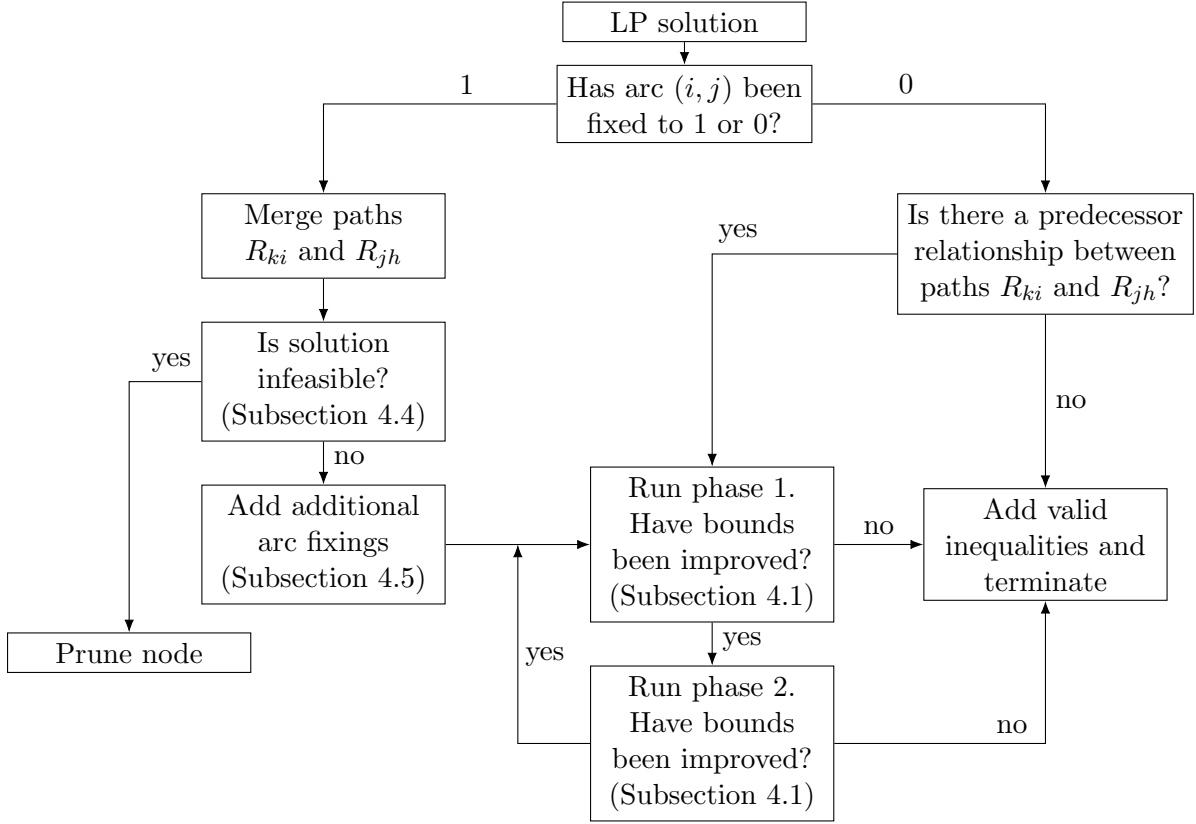
\begin{figure}[htb]%
\centering
\input{flowchart.pgf}%
\caption{Steps in each node of the Branch-and-Cut tree with a fractional LP solution (with the exception of the root node). The valid inequalities include the ones from the literature as well as the ones from Subsection \ref{sec:validInequality}.}%
\label{fig:flowchart}%
\end{figure}

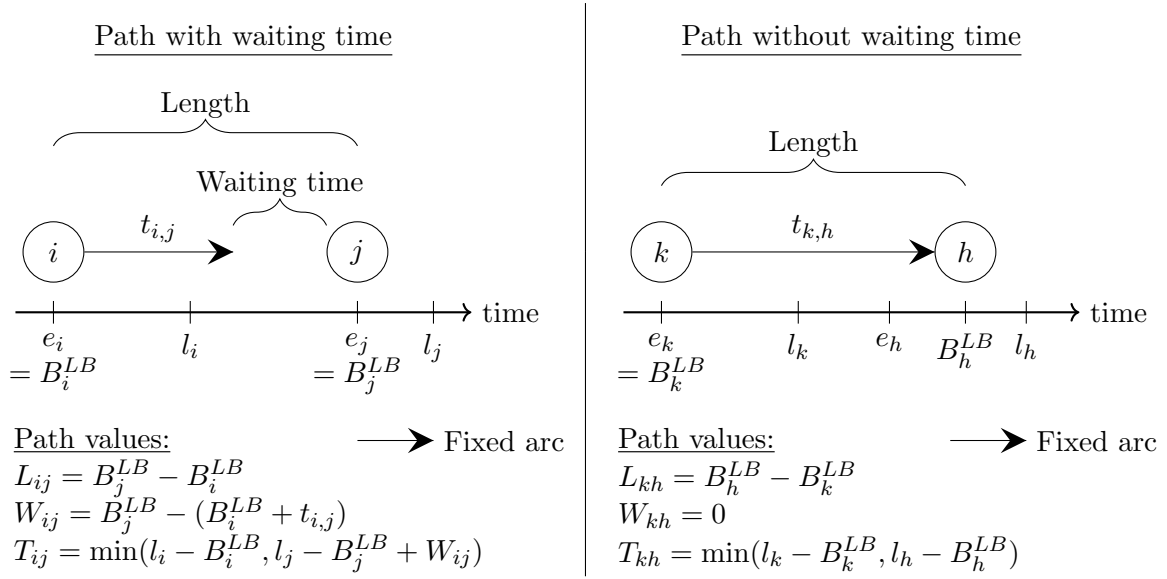
\begin{figure}[htbp]%
\begin{center}
\input{singlePath.pgf}	
\end{center}
\caption{Two fixed paths with two vertices each. Due to the positive waiting time for the path on the left, if vertex $i$ is pushed back, it will first reduce the waiting time before vertex $j$ is affected. Therefore, the forward time slack $T_{ij}$ includes the waiting time in the second term.}%
\label{fig:singlePath}%
\end{figure}

In the following, we explain our new algorithm based on fixed variables called ``fixed path procedure'' as well as give a general outline over our Branch-and-Cut algorithm. Whenever the LP relaxation in a node of the Branch-and-Cut tree returns an integer-feasible solution, Constraints \eqref{eq:SEC} are checked as described in Subsection \ref{sec:model}. If the solution is fractional, we first execute our fixed path procedure and then run several separation routines for valid inequalities from the literature. These are the subtour elimination constraints and generalized order constraints \citep{Cor06}, strengthened capacity constraints and fork constraints \citep{RCL07}, and reachability constraints \citep{L06}. Our implementations for the separation algorithms follow the descriptions of \citet{Cor06} and \citet{RCL07}. For the fork constraints, we consider paths with a maximum size of six locations and for the reachability constraints we consider incompatible groups of up to two locations. An overview about the general procedure in a node is visualized in Figure \ref{fig:flowchart}. The only exception is the root node where we only run through the valid inequalities and do not execute our fixed path procedure, as it requires fixed variables. \\
\ \\
In the remainder of the section, we explain our fixed path procedure: It is a new lower bound improvement procedure for the $B_i$ variables and therefore also for the objective value since the objective is directly determined by the $B_i$ variables. The procedure can also set certain arc variables $X_{ij}$ to 0 and declare fractional solutions infeasible if, based on their current arc fixings, it is not possible to create a feasible integer solution from them. It requires arc variables to be fixed and is therefore most effective in lower nodes of the search tree when many variable fixings have taken place. \\
\ \\
The general idea is based on the precedence implied by the pairing of pickup and delivery vertices as well as the time windows. If, for example, the arc $(i,j+n)$ is fixed to 1 in a node of the Branch-and-Cut tree, we immediately know that vertex $j$ has to be visited prior to vertex $i$, vertex $i+n$ has to be visited at some point after $j+n$, and all of them have to belong to the same tour. This example would lead to the valid inequalities $B_j+t_{ji}\leq B_i$ and $B_{j+n}+t_{j+n,i+n}\leq B_{i+n}$.\\
\ \\
If we are in a node of the Branch-and-Cut tree, let $\bar{X}^0$ be the set of all arcs $(i,j)$ for which the corresponding variable $X_{ij}$ is fixed to 0. Let $R_{ij}$ be the fixed path from $i$ to $j$, i.e. the path where every arc traversed in the path has been fixed to 1. Let also $B_i^{LB}$ be a locally valid lower bound for the variable $B_i$. For every node in the Branch-and-Cut tree we store information about all fixed paths. For every path $R_{ij}$, we store the following information: the duration of the fixed path $L_{ij}$, the cumulated waiting time of the vehicle in the path $W_{ij}$, the forward time slack $T_{ij}$ (the maximum amount of time that we can push back the path without violating a time window, see \citet{S92}, \citet{GD19}) across all vertices in the path, the sequence of all vertices in the path $V_{ij}$ and all direct and indirect predecessor paths $\Pi_{ij}$. A direct predecessor path is a path with at least one pickup vertex to which the corresponding delivery vertex is in path $R_{ij}$. An indirect predecessor path is either a direct or indirect predecessor to a predecessor path. We also save the starting bound, i.e. the earliest time that we can visit the first vertex $i$, which is given by $B_i^{LB}$. We do not explicitly save the time at which we leave the path, i.e. the lower bound for vertex $j$, since it can be calculated in constant time as $B_j^{LB} = B_i^{LB} + L_{ij}$. An example for two fixed paths with two vertices is shown in Figure \ref{fig:singlePath}.\\
\ \\
In Subsection \ref{sec:generalDesc}, we outline the general algorithm. Subsection \ref{sec:updatingPath} details how a path can be updated after we increase the starting bound and Subsection \ref{sec:detour} explains how we can improve the bounding with a shortest path algorithm. In Subsection \ref{sec:identifyInf}, we examine how some of the LP solutions can be pruned before we give examples for fixing additional arc variables in Subsection \ref{sec:fixVar}. Afterwards, we present additional valid inequalities in Subsection \ref{sec:validInequality}. 
\subsection{General procedure}\label{sec:generalDesc}
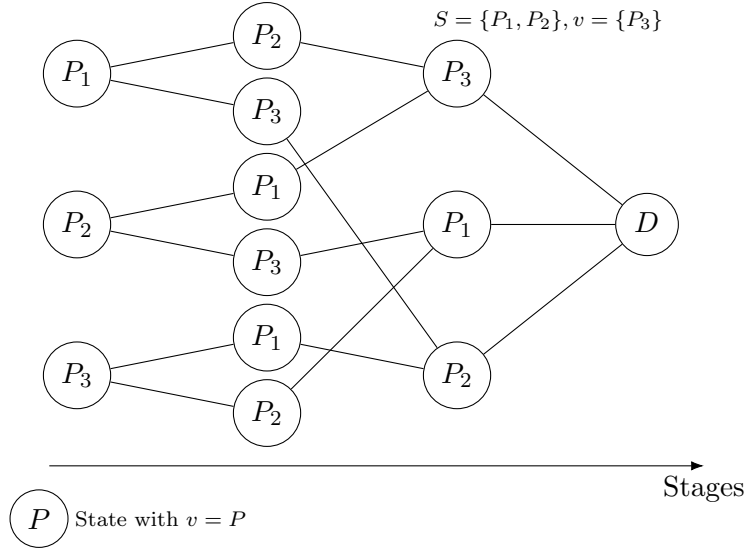
\begin{figure}[htbp]%
\begin{center}
\input{pathResequence.pgf}%
\end{center}
\caption{Example for the dynamic program in the path resequencing. The delivery path (successor path) $D$ is assumed to contain the vertices $\{D_1,D_2,D_3\}$ and has the predecessor paths $S'=\{P_1,P_2, P_3\}$. In this example, all sequences are assumed to be feasible.}%
\label{fig:pathResequencing}%
\end{figure}

In the root node of our Branch-and-Cut tree, there are $2n$ singleton paths $R_{ii}$, one for each location $i\in P\cup D$, with $B_i^{LB} = e_i$, $L_{ii}=0, W_{ii}=0, T_{ii} = l_i-e_i, V_{ii} = \{i\} \mbox{ and either }$
\begin{align*}
&\Pi_{ii}=\{\}\quad &&\mbox{ if } i\in P \mbox{ or }\\
&\Pi_{ii}=\{R_{i-n,i-n}\}\quad &&\mbox{ if } i\in D.
\end{align*}
Whenever an arc $(i,j)$ is fixed to 1 in a node, the two corresponding paths $R_{ki}$ and $R_{jh}$ are merged. Since we know the duration and waiting time for both paths, we can calculate the new lower bound for the end vertex $h$ in constant time. The merging is detailed in Algorithm \ref{code:merging}. If the value $\Delta$ in the pseudo code is positive, we have to wait $\Delta$ units of time when arriving at vertex $j$. Therefore, the time at which we arrive at vertex $h$ does not change. If $\Delta$ is negative, vertex $j$ is pushed back which means vertex $h$ will be pushed back if the waiting time in path $R_{jh}$ is smaller than the amount of push back. If there are other paths with either $R_{ki}$ or $R_{jh}$ as a predecessor, their sets of direct and indirect predecessors have to be updated as well. \\
\begin{algorithm}[htbp]
\caption{Merging of two paths.}
\label{code:merging}
\begin{algorithmic}
\Function{MergePaths}{$R_{ki}$, $R_{jh}$}
\State $\Delta \gets B_j^{LB} - (B_i^{LB} + t_{ij})$\Comment{if $\Delta< 0$: vertex $j$ is pushed back}
\State $R_{kh}\gets \mbox{new path}$
\State $B_h^{LB}\gets B_j^{LB} + L_{jh} - \min(\Delta + W_{jh}, 0)$
\State $L_{kh}\gets B_h^{LB}-B_k^{LB}$
\State $W_{kh}\gets W_{ki} + \max(W_{jh} + \Delta, 0)$
\State $T_{kh}\gets \min(T_{ki}, T_{jh} + W_{ki} + \Delta)$
\State $V_{kh}\gets (V_{ki}, V_{jh})$\Comment{Concatenate the two sequences}
\State $\Pi_{kh}\gets (\Pi_{ki}\cup \Pi_{jh})\backslash \{R_{ki}\}$ \Comment{Remove $R_{ki}$ in case it was predecessor of $R_{jh}$}
\State \Return $R_{kh}$
\EndFunction
\end{algorithmic}
\end{algorithm}
\ \\
After the merging, we run the following two phases: In phase 1, we evaluate the successor paths of the new merged path $R_{kh}$, i.e. all paths $R_{lm}$ with $R_{kh}\in \Pi_{lm}$. When examining a successor path $R_{lm}$, we can potentially improve the lower bound of its starting vertex by setting $B_l^{LB} = \max(B_h^{LB} + tt_{hl}, B_l^{LB})$, where $tt_{hl}$ is the duration of the fastest connection to reach vertex $l$ from vertex $h$. In general, $tt_{hl} = t_{hl}$ holds but if the arc $(h,l)$ is fixed to 0, we can potentially improve this bound by using the shortest detour which is explained in more detail in Subsection \ref{sec:detour}. If we manage to increase the starting bound of path $R_{lm}$, we can check in constant time whether the end vertex's bound $B_m^{LB}$ has been increased as well (see Subsection \ref{sec:updatingPath}). If so, all successors of path $R_{lm}$ can be examined in the next iteration. Once we reach an iteration in which no more successor paths are left to examine, we start phase 2. An example of the bound improvement procedure in phase 1 is given in Figure \ref{fig:phaseOne}. In the following figures, $P_i$ and $D_i$ refer to the pickup and delivery vertex of request $i$, respectively. The path $R_{P_1,P_2}$ in the figure has just been merged so that all of its successors are evaluated. Should $P_3$'s arrival time bound be increased, the successor path containing $D_3$ can then be checked.\\
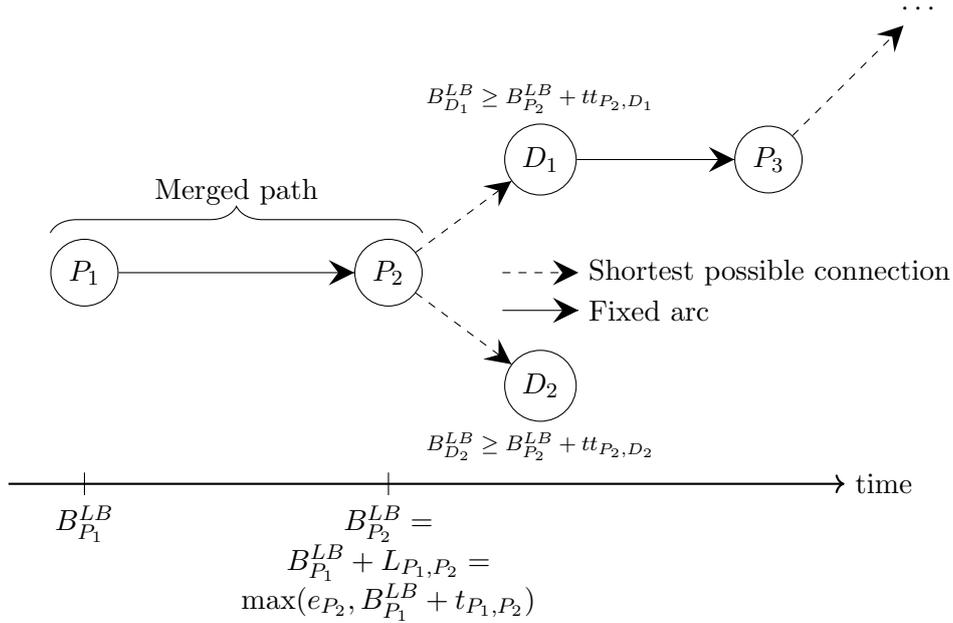
\begin{figure}[htbp]%
\begin{center}
\input{phaseOne.pgf}
\end{center}
\caption{Example for the bound improvement after merging two paths. $P_1$ and $P_2$ were merged, so the connections to all of their successors are evaluated.}%
\label{fig:phaseOne}%
\end{figure}
\ \\
In phase 2, we run a path resequencing procedure which examines all paths that have at least two predecessor paths. We can then evaluate the duration of all permutations of the predecessor paths (both direct and indirect) by also taking into account the connections between the paths $tt_{ij}$. The minimum permutation is the one that leads to the lowest starting bound for the successor path which is a valid lower bound for the successor path's starting time. We find this minimum permutation by running a dynamic program based on the dynamic programming algorithm for the TSP (\citet{Bel62}, \citet{HK62}). Let $S'$ be the set of all paths in the sequence. Then, a state in the dynamic program is represented by the three dimensions $(S,v,t)$ where $S\subseteq S'$ is the set of paths that are already visited (in any order) without $v$, $v\in S'$ is the last visited path (which also means we know the current location) and $t$ is the current time. A state $(S,v,t)$ will always dominate a state $(S,v,t')$ if $t<t'$. This ensures that for every combination of $S$ and $v$ we have at most one state. Additionally, when creating new states, we can take precedence relationships and time windows into account which further reduces the number of sequences that need to be evaluated. If, for example, state $(S,v,t)$ is expanded by visiting a new path $R_{ij}\in S'\backslash(S\cup \{v\})$,  the expansion is only feasible if $R_{ij}$'s time windows are not violated and if $R_{ij}$ is not a successor to any path that is part of the set $S'$ but not yet visited, i.e. $R_{ij}$ cannot be successor to any state in $S'\backslash(S\cup \{v\})$ since it would violate the precedence relationship. An example for this dynamic program can be seen in Figure \ref{fig:pathResequencing}. In each stage of the dynamic program, the states are expanded by adding a new path which is shown in the respective stage. \\
\ \\
Every successor path in the path resequencing that has its ending bound increased is added to the list of paths whose successors should be evaluated and we resume with phase 1. If we finish the path resequencing and no paths have been added, we end our lower bound procedure for this node and return to the normal Branch-and-Cut algorithm. \\
\ \\
If, in the current node of the Branch-and-Cut tree, no arc has been fixed to 1 but instead to 0, we cannot merge any paths. Instead, if arc $(i,j)$ has been fixed to 0, we check whether the path $R_{ki}$ is a predecessor of path $R_{jh}$. If so, we can now utilize the shortest detour $tt_{ij}$ (see Subsection \ref{sec:detour}) to possibly improve the bound for path $R_{jh}$. If we are able to increase the bound for vertex $h$, we can then run our normal 2-phase improvement procedure by checking all successor paths of path $R_{jh}$. \\ 
\ \\
After finishing the procedure in a node of the Branch-and-Cut tree, the improved lower bounds $B_i^{LB}$ are added as local inequalities in the form $B_i\geq B_i^{LB}$. We only add an inequality for the first vertex in each path since Constraints \eqref{eq:time} will propagate the increase along the entire fixed path. When adding such an inequality after an arc was fixed to 0 between two singleton paths, i.e. $R_{ii}$, $R_{i+n,i+n}$ and $(i,i+n) \in \bar{X}^0$ with $i\in P$, care needs to be taken as customer $i$ might still be rejected later. As a result, the inequality for the delivery vertex in such a case takes the form $B_{i+n}\geq B_{i+n}^{LB}(1-Y_i)$. This way, a rejected customer will always have a detour of 0 and will never violate a time window. Customers that are part of a fixed path with at least one other location can never be rejected in that branch since at least one ingoing or outgoing arc of their location is fixed to 1 (compare Equations \eqref{eq:degOut} and \eqref{eq:degIn}). When creating the child nodes for the current node of the Branch-and-Cut tree, we pass all the information about the fixed paths along to the child.

\subsection{Updating a path after increasing the starting bound}\label{sec:updatingPath}
Whenever we set a new starting bound for a path $R_{ij}$, the following values for the path are updated: Let $B^{\mbox{new}}_{i}$ be the new starting bound with $B^{\mbox{new}}_{i}>B_i^{LB}$ and $\Delta=B^{\mbox{new}}_{i}-B_i^{LB}$. Then, $B_i^{LB}\coloneqq B^{\mbox{new}}_{i}$, $L_{ij}\coloneqq L_{ij} - \min(W_{ij}, \Delta)$, $W_{ij}\coloneqq \max(W_{ij}-\Delta,0)$, and $T_{ij}\coloneqq T_{ij}-\Delta$. As $B_j^{LB}$ is not stored explicitly, it does not need to be updated. Should $T_{ij}<0$, the solution is infeasible, as at least one time window in the path is violated. In such a case, we can mark the current node as infeasible and return to our Branch-and-Cut tree.

\subsection{Finding shortest detours between paths}\label{sec:detour}
As mentioned previously, when examining the shortest connection $tt_{ij}$ between two paths $R_{ki}$ and $R_{jh}$, we can improve the bounding if $(i,j) \in \bar{X}^0$, as we know that a detour through at least one other path has to be used instead. We try to find the minimum arrival time at vertex $j$ by solving a modified shortest path problem. In this shortest path problem, the paths can be understood as vertices and the arcs between the vertices correspond to the arcs between the paths. The final result is not necessarily a feasible shortest detour but guaranteed to be a lower bound for the shortest detour. A more detailed description of the algorithm can be found in Appendix \ref{sec:app_Detour}.

\subsection{Identifying infeasible solutions ahead of time}\label{sec:identifyInf}
Whenever we fix a new arc and therefore merge two paths, we run a short algorithm to detect any circular precedence relations that might arise from the merge. We realize this by checking whether any of the direct and indirect predecessor paths of the new merged path have the new merged path as their own predecessor. In such a case, we can mark the current node as infeasible. We can also mark a solution as infeasible if the new path is connected to the start depot but has predecessor paths. This works analogously for paths connected to the end depot and potential successors.\\

\subsection{Fixing additional arc variables}\label{sec:fixVar}
In addition to improving the bounds for the $B_i$ variables, we also fix certain arc variables to 0 based on the information in the paths. Whenever two paths are merged to create one new path $R_{ij}$, we can immediately set the variable $X_{ji}$ to 0, as it would lead to a cycle. The same applies to any arcs from $j$ to predecessor paths $\Pi_{ij}$ and arcs from any successor paths to $i$ which would result in circular relationships. Additionally, for any path $R_{ij}$ with at least one predecessor (successor) path, the arc $X_{0i}$ ($X_{j,2n+1}$) has to be 0. Furthermore, we can fix arc variables whenever there is a path that either starts at the depot (i.e. there is at least one path $R_{0i}$ for any $i\in P\cup D$) or ends at the depot (i.e. at least one path $R_{i,2n+1}$ for any $i\in P\cup D$). In the latter case, we can evaluate all arcs from direct and indirect predecessor paths $\Pi_{jh}$ of another path $R_{jh}$. All arcs 
\begin{align*}
X_{ki}\quad&\forall R_{lk}\in \Pi_{jh} 
\end{align*}
can be set to 0, as they would always lead to a violation of the precedence constraints in the tours. The procedure works analogously for a path connected to the start depot and successors of other paths. \\
\ \\
Any arcs connecting two paths that would lead to a violation of a time window in the second path can also be set to 0. The same goes for any arcs connecting two paths where the passengers would exceed the capacity of the vehicles. The capacity for the path sequence $\{R_{ij}, R_{kh}\}$ can be checked by running through both paths' vertex sequences, i.e. $V_{ij}$ and $V_{kh}$, and adding up the passengers that enter or exit in each vertex. Should the maximum number of passengers ever surpass the vehicle capacity, the two paths cannot be visited straight after one another and $X_{jk}$ can be set to 0. If a delivery vertex is found in the second path $R_{kh}$ and the corresponding pickup was not visited in either path, the passengers for that request must have been in the vehicle the entire time (i.e. they entered before $R_{ij}$), so the previously found maximum can be increased correspondingly.\\
\ \\
If the merged path $R_{ij}$ has a successor path $R_{kh}$, we also check all sequences including one intermediate path, i.e. the sequence $\{R_{ij}, R_{lm}, R_{kh}\}$ for every path $R_{lm}$. If any of the time windows in path $R_{lm}$ or $R_{kh}$ are violated, the arc variables $X_{jl}$ and $X_{mk}$ can be set to 0. This does not apply to violations of the capacity, as the number of passengers could be reduced by including additional paths with delivery vertices in between. \\
\ \\
Like \citet{Cor06} we identify all pairs of incompatible customers in the preprocessing, i.e. customers who can, because of time windows or passenger numbers, never be feasibly visited in the same tour together. Whenever a path $R_{ki}$ is merged, we can determine all customers now incompatible with the path based on the locations in the path $V_{ki}$. Then, all arcs to and from other paths containing one of the incompatible customers can be set to 0.\\

\subsection{Additional valid inequalities}\label{sec:validInequality}
Before adding the current lower bounds for the $B$ variables at the end of our procedure, we check the violation of some more constraints to improve the big M formulation of Constraints (\ref{eq:time}), as big M constraints often lead to very poor LP relaxations and should be formulated as tight as possible (\citet{CRT90}, \citet{CF06}). For every path $R_{jh}$, the following inequality is valid:
\begin{equation}
\sum_{i\in P\cup D} (B_{i}^{LB} + t_{ij})\cdot X_{ij} \leq B_{j}\label{eq:addB1}
\end{equation}
The inequality can be made less dense by only considering the vertices $i$ that are the end of a path, i.e. all paths $R_{ki}$. Furthermore, when examining a pair of paths $R_{ki}$ and $R_{i+n,m}$, the following inequality is valid:
\begin{equation}
\sum_{R_{jh}}(B_i^{LB} + t_{ij} + L'_{jh} + tt_{h,i+n})\cdot X_{ij}\leq B_{i+n} \label{eq:addB2}
\end{equation}
where $L'_{jh}$ is the updated length of path $R_{jh}$ that might arise from a later arrival at $j$ (compare Subsection \ref{sec:updatingPath}). Effectively, we add to the earliest possible departure from vertex $i$ the travel time to reach vertex $j$, then the duration of the path to get to vertex $h$ and then the travel time (including possible detours) when going from vertex $h$ to $i+n$ and weight this sum with the current flow on the arc $(i,j)$. These inequalities are only checked for violations if at least one $B_i^{LB}$ was improved.

\section{Computational study}\label{sec:comp}
The algorithms were coded in C++ and compiled with Visual Studio 2017 (version 15.9.8). We used the CPLEX API with CPLEX version 12.9 and the algorithms are implemented via callbacks. Unless otherwise noted, all CPLEX parameters are left at their default values. The tests were executed on a machine with an AMD Ryzen Threadripper 3990X with 2.9GHz and 256GB RAM. Each instance was limited to a single thread. The starting solution for each run is a trivial solution in which every customer is rejected.\\
\ \\
All used instances can be found under the following link: \url{http://doi.org/10.25592/uhhfdm.10389}.

\subsection{Instances}
The utilized instances were randomly generated but are based on a real-life ridepooling service in Hamburg, Germany. The customer locations for each instance were drawn randomly from the stops of the ridepooling provider. The resulting distance matrix is not symmetrical but fulfills the triangle inequality. On average, the direct travel time between a pickup and its corresponding delivery location is between eight and nine minutes. For each customer $i\in P$ the start of the pickup time window $e_i$ was drawn such that the difference between the pickup time window starts is exponentially distributed with $\lambda=0.005$ which means the time window beginnings differ by on average 200 seconds. Additionally, $e_i\geq\max_{i\in P}(t_{0i})$, i.e. each pickup time window can be reached on time from the depot. The end of the pickup time window $l_i$ was chosen such that the time window has a length of 5, 10 or 15 minutes (equally weighted). The delivery time window starts at $l_{i+n}=e_i+t_{i,i+n}$, i.e. we never have to wait at a delivery vertex, and the time window end follows as $l_{i+n}=l_{i} + t_{i,i+n}\cdot u$ where $u\sim U(\alpha-0.1,\alpha+0.1)$ and $\alpha$ is a factor that increases the length of the delivery time window to allow for more combinations of different customers in a tour. For example, if $\alpha = 1.3$, it means that if a customer's pickup location is left at the latest possible time $l_{i}$, the vehicle does not immediately have to visit the delivery location but instead has a buffer of $20-40\%$ of the direct travel time. \\
\ \\
The passenger capacity for each vehicle was set to six, which is also common in practice. The number of passengers per request $q_i$ was drawn according to an exponential distribution with $\lambda=0.9$ and rounded up. If the resulting value was greater than the vehicle capacity of six, the value was discarded and redrawn.\\
\ \\
For our computational study we generated instances for different combinations of $n$ and $K$. The values for $n$ ranged between 10 and 120 and for $K$ between 1 and 10. For each value of $n$ we chose five values of $K$. The combinations can be seen in Table \ref{tab:instances}. Additionally, we varied the time window length factor $\alpha$ between the values $\{1.1,1.3,1.5\}$. For each combination $\{n,K,\alpha\}$ we generated 20 instances which results in a testbed of 4500 instances in total. 

\begin{table}[htbp]%
\centering
\caption{The combinations of $n$ and $K$ for which we generated instances.}
\label{tab:instances}
\begin{tabular}{lcccccccccc}
\toprule
$n \backslash K$ & 1 & 2 & 3 & 4 & 5 & 6 & 7 & 8 & 9 & 10\\
\midrule
10&$\bullet$ &$\bullet$ &$\bullet$ &$\bullet$ &$\bullet$ &&&&&\\
15&$\bullet$ &$\bullet$ &$\bullet$ &$\bullet$ &$\bullet$ &&&&&\\
20&&$\bullet$ &$\bullet$ &$\bullet$ &$\bullet$ &$\bullet$ &&&&\\
25&&$\bullet$ &$\bullet$ &$\bullet$ &$\bullet$ &$\bullet$ &&&&\\
30&&&$\bullet$ &$\bullet$ &$\bullet$ &$\bullet$ &$\bullet$ &&&\\
35&&&$\bullet$ &$\bullet$ &$\bullet$ &$\bullet$ &$\bullet$ &&&\\
40&&&&$\bullet$ &$\bullet$ &$\bullet$ &$\bullet$ &$\bullet$ &&\\
45&&&&$\bullet$ &$\bullet$ &$\bullet$ &$\bullet$ &$\bullet$ &&\\
50&&&&&$\bullet$ &$\bullet$ &$\bullet$ &$\bullet$ &$\bullet$ &\\
60&&&&&&$\bullet$ &$\bullet$ &$\bullet$ &$\bullet$ &$\bullet$ \\
70&&&&&&$\bullet$ &$\bullet$ &$\bullet$ &$\bullet$ &$\bullet$ \\
80&&&&&&$\bullet$ &$\bullet$ &$\bullet$ &$\bullet$ &$\bullet$\\
90&&&&&&$\bullet$ &$\bullet$ &$\bullet$ &$\bullet$ &$\bullet$\\
100&&&&&&$\bullet$ &$\bullet$ &$\bullet$ &$\bullet$ &$\bullet$\\
120&&&&&&$\bullet$ &$\bullet$ &$\bullet$ &$\bullet$ &$\bullet$\\\bottomrule
\end{tabular}
\end{table}
\ \\
Each run was limited to one hour of computational time.

\subsection{Component tests}
Since our procedure consists of many components, we ran a preliminary series of tests where single components are turned off to see whether their inclusion is actually worthwhile. For this test we created six configurations: The first one is a benchmark configuration where every component is turned on. In the second configuration, we turned off our entire fixed path procedure and only separated the inequalities from the literature in a node of the Branch-and-Cut tree. In the remaining configurations, we turned off one component each: The additional inequalities for the $B_i$ variables (Inequalities \eqref{eq:addB1} and \eqref{eq:addB2}), the path resequencing (phase 2 in Subsection \ref{sec:generalDesc}), the calculation of the shortest detours (Subsection \ref{sec:detour}), and all of our own preprocessing fixings and inequalities (Subsection \ref{sec:preproc}). All configurations still include the described steps from \citet{Cor06} and \citet{RCL07}.\\
\ \\
The instances for these tests are generated in the same way as our other instances but we used a different random seed in the creation so as to keep the sets separate. These new instances include $n\in \{90,100\}$ with the five different vehicle levels, the three time window lengths and 20 instances per combination, i.e. 600 instances in total.\\
\ \\
The average final gap for each configuration is shown in Table \ref{tab:components}. We also tested whether the gaps were statistically significantly different from the benchmark configuration. All of the tested components have a positive influence since their removal significantly worsens the average gap. The biggest impact stems from the removal of the preprocessing steps. The configuration ``fixed paths'' effectively turns off all the three components ($B_i$ inequalities, path resequencing, and shortest detour). Therefore, it is unsurprising that it performs worse than the configurations where they are turned off individually. Of those three components, the $B_i$ inequalities appear to improve the gaps the most.

\begin{table}[htb]
\centering
\caption[]{The average MIP gap for each configuration and an indication whether the gaps are statistically different from the benchmark configuration.\\ $^*$: $p\leq0.05$, $^{**}$: $p\leq 0.001$}
\label{tab:components}
\begin{tabular}{lr}
  \toprule
Configuration & avg. gap \\ 
  \midrule
benchmark & 0.0612\phantom{$^{**}$} \\ 
  fixed paths & 0.0900$^{**}$ \\ 
  $B_i$ inequalities & 0.0790$^{**}$ \\ 
  path resequencing & 0.0713$^{**}$\\ 
  shortest detour & 0.0655$^*$\phantom{$^{*}$} \\ 
  preprocessing & 0.1813$^{**}$ \\ 
   \bottomrule
\end{tabular}
\end{table}

\subsection{Results}
For our final calculations, we considered two configurations: \textit{enhanced} is the configuration with all the described preprocessing and fixed path inequalities enabled while \textit{basic} is a benchmark configuration that does not utilize our own developed preprocessing inequalities nor any of the algorithms of the fixed path procedure. However, it still uses the other improvements we added from the literature (including the valid inequalities), so that we can accurately measure the effect of our methods. \\
\ \\
An overview of the results of the two configurations for the different values of $n$ is presented in Table \ref{tab:configs_all}. For smaller instances ($n=10, 15$) both configurations can find all optimal solutions, but for larger values of $n$ fewer runs can prove optimality in the given time limit. \textit{enhanced} can clearly outperform \textit{basic} in the number of proven optimal solutions found as well as upper and lower bounds. Only for four values of $n$ (35, 45, 50, and 90) the average objective value of \textit{basic} is better. These differences are largely caused by a higher number of rejected customers in some of the non-optimal solutions which raises the average significantly. Altogether, \textit{basic} was able to find a better upper bound in only 110 instances and a better lower bound in 13 instances. On the other hand, \textit{enhanced} provided better upper bounds in 758 instances and better lower bounds in 1646 instances. The number of instances solved in the root node is higher for \textit{enhanced} as well. This indicates that our preprocessing inequalities are helpful in bounding the problem, since in these instances our fixed path procedure was never executed.

\begin{table}[htbp]
\centering
\caption[]{Comparison between the two configurations over all instances.\\ opt: number of instances solved to optimality, root: number of instances solved in the root node, \textit{z}: average objective value, LB: average lower bound, time: average time used in seconds.}
\label{tab:configs_all}
\begin{tabular}{lrrrrrrrrrr}
  \toprule
	&\multicolumn{5}{c}{\textit{enhanced}} & \multicolumn{5}{c}{\textit{basic}}\\\cmidrule(lr{.75em}){2-6}\cmidrule(lr{.75em}){7-11}
\multicolumn{1}{c}{$n$} & \multicolumn{1}{c}{opt}  & \multicolumn{1}{c}{root} & \multicolumn{1}{c}{\textit{z}} & \multicolumn{1}{c}{LB} & \multicolumn{1}{c}{time} &\multicolumn{1}{c}{opt} &\multicolumn{1}{c}{root} & \multicolumn{1}{c}{\textit{z}} & \multicolumn{1}{c}{LB} & \multicolumn{1}{c}{time}  \\ 
  \midrule
10 & 300 & 160 & 85.62 & 85.62 & 0.19 & 300 & 107 & 85.62 & 85.62 & 0.28 \\ 
  15 & 300 &  66 & 210.56 & 210.56 & 4.16 & 300 &  38 & 210.56 & 210.56 & 12.45 \\ 
  20 & 298 &  61 & 157.30 & 156.83 & 77.03 & 292 &  35 & 157.35 & 156.16 & 187.27 \\ 
  25 & 276 &  44 & 277.95 & 266.30 & 482.60 & 256 &  20 & 278.36 & 262.32 & 772.75 \\ 
  30 & 270 &  59 & 135.18 & 113.45 & 556.84 & 236 &  38 & 135.76 & 107.97 & 1013.88 \\ 
  35 & 245 &  49 & 204.18 & 153.82 & 859.85 & 197 &  22 & 202.59 & 149.30 & 1479.82 \\ 
  40 & 247 &  88 & 70.11 & 37.31 & 798.10 & 177 &  51 & 75.69 & 33.08 & 1602.47 \\ 
  45 & 228 &  78 & 76.03 & 48.70 & 1031.83 & 149 &  52 & 73.11 & 43.99 & 1982.81 \\ 
  50 & 235 & 117 & 35.34 & 11.76 & 889.77 & 155 &  75 & 34.01 & 10.50 & 1822.08 \\ 
  60 & 248 & 163 & 7.03 & 3.32 & 689.29 & 177 & 111 & 9.03 & 2.54 & 1535.87 \\ 
  70 & 230 & 137 & 3.35 & 2.90 & 905.31 & 151 & 103 & 162.32 & 2.26 & 1826.56 \\ 
  80 & 215 & 130 & 13.27 & 3.49 & 1089.59 & 128 &  87 & 145.16 & 2.64 & 2115.78 \\ 
  90 & 201 & 126 & 132.61 & 3.65 & 1274.97 & 115 &  77 & 30.48 & 2.67 & 2250.21 \\ 
  100 & 182 & 110 & 16.41 & 4.01 & 1480.82 &  96 &  67 & 628.65 & 2.90 & 2472.20 \\ 
  120 & 163 &  85 & 68.91 & 4.87 & 1711.83 &  71 &  46 & 2947.87 & 3.44 & 2778.83 \\ 
   \bottomrule
\end{tabular}
\end{table}

\ \\
Looking at the average computation time, \textit{enhanced} is faster. However, when comparing the computation times, looking at the average of all instances is slightly unfair towards \textit{basic}, since we already showed that the configuration hits the time limit more often. Therefore, Table \ref{tab:optimal} contains the average computation times of all instances that were solved to optimality by both configurations. It is evident that \textit{enhanced} can prove optimality faster than \textit{basic} even when only looking at this subset of instances. Focusing only on \textit{enhanced}, Table \ref{tab:times} shows the computation time spent on our fixed path procedure during the entire Branch-and-Cut algorithm and the total time needed. As a whole, the time spent for our procedure is negligible compared to the time spent on the entire problem. This holds even when taking into account the maximum time spent for each $n$ which is unbiased by the fast completion time of the instances that are solved in or close to the root node. Our procedure is therefore  lightweight and provides benefit at little cost.

\begin{table}[htbp]
\centering
\caption{The computational times of all instances for which both configurations found the optimal solution.}
\label{tab:optimal}
\begin{tabular}{lrrr}
  \toprule
\multicolumn{1}{c}{$n$} & \multicolumn{1}{c}{\# instances} & \multicolumn{1}{c}{\textit{enhanced} time (s)}  &\multicolumn{1}{c}{\textit{basic} time (s)} \\ 
  \midrule
10 & 300 & 0.19 & 0.28 \\ 
  15 & 300 & 4.16 & 12.45 \\ 
  20 & 292 & 31.24 & 93.77 \\ 
  25 & 256 & 104.70 & 286.82 \\ 
  30 & 236 & 78.01 & 312.56 \\ 
  35 & 197 & 73.60 & 371.28 \\ 
  40 & 177 & 22.68 & 214.16 \\ 
  45 & 149 & 35.31 & 343.74 \\ 
  50 & 155 & 8.80 & 158.79 \\ 
  60 & 177 & 1.27 & 101.40 \\ 
  70 & 151 & 2.16 & 76.44 \\ 
  80 & 128 & 2.38 & 121.17 \\ 
  90 & 115 & 1.63 & 78.40 \\ 
  100 &  96 & 2.05 & 75.09 \\ 
  120 &  71 & 1.54 & 128.72 \\ 
   \bottomrule
\end{tabular}
\end{table}

\begin{table}[htb]
\centering
\caption{The computational times of just the fixed path procedure (FPP) and of the entire run in comparison as well as the number of calls to the fixed path procedure.}
\label{tab:times}
\begin{tabular}{lrrrr}
  \toprule
 &\multicolumn{3}{c}{time (s)} &\\\cmidrule(lr{.75em}){2-4}
 	$n$&\multicolumn{1}{c}{FPP (mean)} & \multicolumn{1}{c}{FPP (max)} & \multicolumn{1}{c}{total run (mean)} & \multicolumn{1}{c}{\# calls of FPP} \\ 
  \midrule
10 & 0.00 & 0.00 & 0.19 & 12.69 \\ 
  15 & 0.00 & 0.00 & 4.16 & 355.08 \\ 
  20 & 0.00 & 0.00 & 77.03 & 4172.06 \\ 
  25 & 0.00 & 0.02 & 482.60 & 16619.98 \\ 
  30 & 0.00 & 0.00 & 556.84 & 15553.73 \\ 
  35 & 0.00 & 0.01 & 859.85 & 20046.49 \\ 
  40 & 0.00 & 0.01 & 798.10 & 17343.69 \\ 
  45 & 0.00 & 0.04 & 1031.83 & 26901.92 \\ 
  50 & 0.00 & 0.06 & 889.77 & 25872.36 \\ 
  60 & 0.02 & 1.14 & 689.29 & 21197.50 \\ 
  70 & 0.09 & 1.87 & 905.31 & 18188.45 \\ 
  80 & 0.41 & 12.30 & 1089.59 & 16818.86 \\ 
  90 & 1.02 & 10.71 & 1274.97 & 14628.61 \\ 
  100 & 1.46 & 15.10 & 1480.82 & 14771.69 \\ 
  120 & 1.45 & 15.61 & 1711.83 & 10253.86 \\ 
   \bottomrule
\end{tabular}
\end{table}
\ \\
Dividing the instances based on the time window length factor $\alpha$ is another way of measuring their difficulty. In Table \ref{tab:timewindows}, the sum of optimal solutions for the three different settings of $\alpha$ are shown. Our results show that a setting with smaller time windows (i.e. a smaller $\alpha$) makes instances easier to solve, as more parts of the solution space can be disregarded. This effect holds even though the tighter time windows do not lead to infeasibility in our case, since a rejection is still possible. Expectedly, larger time windows generally lead to fewer rejections though this effect seems to hold mostly for the smaller instances where the rejection rates are also higher in general. Since the instances with more requests also use more tours, the rate of rejections declines, as more tours allow for a better distribution of requests (e.g. cluster similar requests together). 

\begin{figure}[htbp]%
\includegraphics[width=\columnwidth]{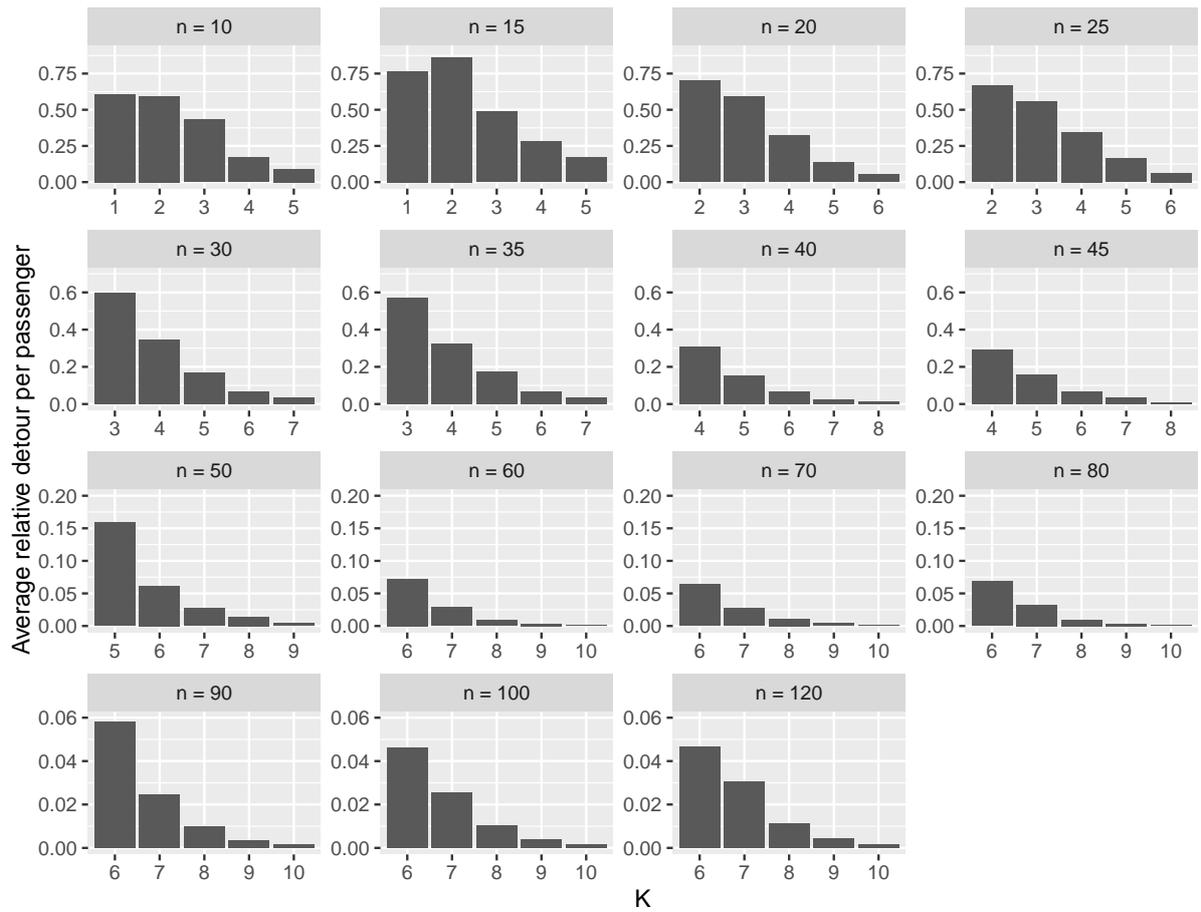}%
\caption{The average relative detour of all non-rejected passengers in the optimal solutions of configuration \textit{enhanced}. Note that the vertical axis scale is constant per row but changes across rows.}%
\label{fig:averageObjective}%
\end{figure}

\begin{table}[htbp]
\centering
\caption[]{The instances split by the time window length factor $\alpha$.\\ opt: number of optimal solutions, rejected: number of instances where the solution contains at least one rejected customer.}
\label{tab:timewindows}
\begin{tabular}{lrrrrrr}
  \toprule
	& \multicolumn{2}{c}{$\alpha = 1.1$} & \multicolumn{2}{c}{$\alpha = 1.3$} & \multicolumn{2}{c}{$\alpha = 1.5$}\\\cmidrule(lr{.75em}){2-3}\cmidrule(lr{.75em}){4-5}\cmidrule(lr{.75em}){6-7}
\multicolumn{1}{c}{$n$} & \multicolumn{1}{c}{opt} & \multicolumn{1}{c}{rejected} & \multicolumn{1}{c}{opt}& \multicolumn{1}{c}{rejected} & \multicolumn{1}{c}{opt} & \multicolumn{1}{c}{rejected} \\ 
  \midrule
10 & 100 &  53 & 100 &  48 & 100 &  45 \\ 
  15 & 100 &  62 & 100 &  58 & 100 &  56 \\ 
  20 & 100 &  49 & 100 &  40 &  98 &  38 \\ 
  25 & 100 &  51 &  95 &  52 &  81 &  47 \\ 
  30 &  93 &  40 &  93 &  28 &  84 &  28 \\ 
  35 &  92 &  39 &  79 &  32 &  74 &  33 \\ 
  40 &  90 &  19 &  84 &  17 &  73 &  13 \\ 
  45 &  84 &  26 &  73 &  13 &  71 &  13 \\ 
  50 &  84 &   8 &  80 &   7 &  71 &   6 \\ 
  60 &  88 &   3 &  79 &   0 &  81 &   2 \\ 
  70 &  83 &   0 &  78 &   0 &  69 &   0 \\ 
  80 &  79 &   2 &  71 &   2 &  65 &   1 \\ 
  90 &  72 &   0 &  71 &   1 &  58 &   0 \\ 
  100 &  63 &   0 &  61 &   2 &  58 &   2 \\ 
  120 &  68 &   1 &  53 &   3 &  42 &   4 \\  
   \bottomrule
\end{tabular}
\end{table}
\ \\
In Figure \ref{fig:averageObjective}, we visualize the development of the average detour per passenger. For that, the first part of the objective value (excluding the penalty factor $\Phi$) is divided by the number of passengers that were served (i.e. not rejected). The figure only includes solutions that are proven optimal, since otherwise a poor quality solution with many rejected customers would unfairly reduce the average detour as higher rejections lead to lower detours for the remaining served customers. As an example, a passenger in one of the instances with $n=10$ and $K=1$ has to endure a detour of roughly 60\% while a passenger in the instances with $n=120$ and $K=6$ of below 6\%. The figure shows that, expectedly, the average detours decrease as the number of vehicles increases. This effect is more pronounced in the instances with more customer requests. This is due to the fact that in the instances with fewer customers, the number of rejections is higher (compare Table \ref{tab:timewindows}). Therefore, an increase in vehicles in instances with small $n$ will primarily increase the number of served requests and only secondarily decrease the average detours. Hence, the size of the vehicle fleet is a crucial parameter, as small variations can have a drastic effect on the average solution quality. By selecting the number of used vehicles, the ridepooling provider can steer the average relative detour of the customers and by this solve the tradeoff between a sufficiently high pooling rate and a sufficiently small detour for the customers and thus be an attractive option for the customers as mentioned in the introduction.

\section{Conclusion}\label{sec:conclusion}
The paper investigates a dial-a-ride problem focusing on the residents of urban areas. These citizens can use a wide variety of transportation modes. Because of this, ridepooling providers have to offer short detours for the customers to be sufficiently attractive for them. However, ridepooling providers need high pooling rates to operate economically efficient. We developed a Branch-and-Cut algorithm to determine tours with as short as possible detours given the number of vehicles while serving as many customers as possible. Our computational results show that ridepooling providers can steer the average relative detour of the customers by selecting the number of used vehicles and therefore solve the tradeoff mentioned above according to their needs. \\
\ \\
We introduced a new technique using information about already fixed paths in the Branch-and-Cut nodes to improve lower bounds on the arrival time variables at the customer locations. As our objective function (\ref{eq:objFunc}) minimizes relative detours using arrival times at the delivery locations, our improved lower bounds directly improve the local objective function's lower bound in the node. We used them further to introduce additional valid inequalities. Moreover, we used the paths to identify infeasible solutions ahead of time. The computational study showed that our procedure performs significantly better than the mixed-integer programming formulation including the standard valid inequalities. \\
\ \\
We used the information that a customer's pickup location has to be visited by the same vehicle and before its delivery location to determine our paths. As a direction for future research our technique using the paths information can be used to improve the solution process of different routing problems with precedence relationships. \\

\section*{Declaration}
\section*{Funding}
This project was supported by the Hamburger Behörde für Wissenschaft, Forschung, Gleich\-stellung und Bezirke (BWFGB; Hamburg authority for science, research, equalization, and districts). No grant number is available.

\section*{Conflict of interest}
The authors declare that they have no conflict of interest.

\section*{Availability of data}
The instances are available under the following link: \url{http://doi.org/10.25592/uhhfdm.10389}.

\bibliographystyle{apa}
\bibliography{DARP_MIP_TW}

\appendix
\section{Finding shortest detours}\label{sec:app_Detour}
Here we briefly describe our algorithm to find the shortest detour between a path $R_{ki}$ and a path $R_{jh}$ if the arc $(i,j)$ has been fixed to 0. In the following, we refer to $R_{ki}$ as the source and $R_{jh}$ as the sink. All paths that are neither predecessors to $R_{ki}$ nor successors to $R_{jh}$ are vertices in this shortest path problem. The arcs between the vertices exist when the respective upper bound of the arc variable that connects the two paths is 1. The length of each path is added to the duration of each outgoing arc, e.g. the arc $(m,j)$ that connects a path $R_{lm}$ to $R_{jh}$ has a duration of $t_{mj} + L_{lm}$. The problem is then solved by a modified Dijkstra Algorithm \citep{D59}. The procedure is shown in Algorithm \ref{code:detour}. In the pseudo code, all path parameters such as the length $L$ and the time slack $T$ are continuously updated when a path's arrival is changed as described in Section \ref{sec:updatingPath}, though all changes to the intermediate paths are local to the function. \\
\ \\
First, all paths that are direct neighbors of the source, i.e. all paths $R_{lm}$ for which $(i,l)\notin \bar{X}^0$, are visited and their respective arrival time is saved. Similar to the Dijkstra Algorithm, each visited path is saved in a priority queue where a lower arrival time constitutes a higher priority. In each iteration, we remove the highest priority path $R_{lm}$ and first evaluate whether there is a direct connection from this path to the sink, i.e. whether $(m,j)\notin \bar{X}^0$ (line \ref{code:detour:IF_Direct}). If so, we can potentially update the minimum arrival at $j$ and then discard the current path from further evaluation since due to the triangle inequality, the direct arc can never be slower than a detour over another path. If the direct arc to the sink is not available, we visit all other neighbour paths of $R_{lm}$ where the time window of the neighbour path is not violated (line \ref{code:detour:neighbor}). These paths are added to the priority queue if their arrival time is now less than before. Checking the time windows (via the time slack $T_{no}$) is required since the original bounds are initialized with $\infty$ and therefore an earlier arrival at the path does not guarantee time window feasibility. Once we reach an iteration where the highest priority path has a higher arrival time than the sink (line \ref{code:detour:abort}) or no more paths are in the queue, we can abort, as no faster path can be found in the network. The final shortest path is not guaranteed to be feasible, as we do not check for passenger capacity constraints nor precedence relationships between the members of the path. However, it is still a valid lower bound for the arrival at the sink.

\begin{algorithm}[htb]
\caption{Finding the shortest detour between two paths.}
\label{code:detour}
\begin{algorithmic}[1]
\Function{FindDetour}{$R_{ki}$, $R_{jh}$}
\State $V\gets\mbox{set of all eligible paths}$
\State $B_v\gets\mbox{array with earliest arrival time per path (initialized with $\infty\ \forall v\in V$)}$
\State $Q\gets\mbox{priority queue (earlier arrival time $\rightarrow$ higher priority)}$
\State $B_{R_{ki}}\gets B_k^{LB}$\label{code:detour:initTime}
\State Add $R_{ki}$ to $Q$
\While{$Q$ not empty}
	\State $R_{lm} \gets \mbox{highest priority path in Q}$
	\If{$B_{R_{lm}} \geq B_{R_{jh}}$}\label{code:detour:abort}
		\State Exit while \Comment{None of the remaining paths can improve $B_{R_{jh}}$}
	\EndIf	
	\If{$(m,j)\notin \bar{X}^0$}\label{code:detour:IF_Direct} 
		\State $B_{R_{jh}} \gets \min(B_{R_{jh}}, B_{R_{lm}} + L_{lm}+ t_{mj})$
	\Else
		\ForAll{paths $R_{no}$ in V}
			\State $B^{\mbox{new}} \gets B_{R_{lm}} + L_{lm} + t_{mn}$\Comment{Arrival at $R_{no}$ after leaving $R_{lm}$}
			\If{$(m,n)\notin \bar{X}^0$ AND $B^{\mbox{new}} < B_{R_{no}}$ AND $B^{\mbox{new}} \leq B^{LB}_n + T_{no}$}\label{code:detour:neighbor}
				\State $B_{R_{no}} \gets \max(B^{\mbox{new}}, e_n)$
				\State add $R_{no}$ to $Q$ with priority $B_{R_{no}}$
			\EndIf
		\EndFor
	\EndIf	
\EndWhile
\State \Return $B_{R_{jh}}$\Comment{Return the earliest possible arrival time at vertex $j$}
\EndFunction
\end{algorithmic}
\end{algorithm}
\ \\
If we call this procedure in the path resequencing in phase 2, we can exclude additional paths from consideration since any path visited in the current sequence as well as predecessors of the paths preceding $R_{ki}$ and successors of the paths after $R_{jh}$ can never be used as a feasible detour which reduces the number of eligible paths. Furthermore, the earliest arrival time at path $R_{ki}$ might be higher than $B_k^{LB}$ if the path was already pushed back in the sequence which allows for an adjustment in the initialization in line \ref{code:detour:initTime} by increasing $B_{R_{ki}}$.\\
\ \\
This shortest path algorithm clearly runs in polynomial time with respect to the number of customers. The number of vertices in the graph is linear in the number of customers in the original problem $n$ since there are at most $2n$ paths. Each path is examined at most once and for each examined path we run through all neighbour paths in the worst case. Finally, the priority queue has a logarithmic worst case complexity to add an item which means our algorithm runs in $O(n^2\log(n))$. This polynomial run time is in contrast to the generally hard shortest path problem with time windows (see \citet{D94}, \citet{IGSD98}). This difference can be explained by the fact that from a dynamic programming point of view, we only need the arrival time to label a vertex and an earlier arrival time at a vertex always dominates a later one. This reasoning is also why we do not consider capacity or precedence constraints between the vertices in the graph, as it would require additional dimensions to solve the problem exactly.

\end{document}

%% file: matching.pgf
\begin{tikzpicture}
\usetikzlibrary{arrows,decorations.markings, decorations.pathreplacing, calc}
\usetikzlibrary{math}
\usetikzlibrary{shapes.misc}
\usetikzlibrary{positioning}

\tikzset{Path/.style = {draw,circle, minimum size=1.4cm}}
\tikzset{Arr/.style = {decoration={markings,mark=at position 1 with
    {\arrow[scale=3,>=stealth]{>}}},postaction={decorate}}}

\definecolor{myBlue}{RGB}{100,143,255}
\definecolor{myPurple}{RGB}{120,94,240}
\definecolor{myMagenta}{RGB}{220,38,127}
\definecolor{myOrange}{RGB}{254,97,0}

\def\baseHeight{2};
\def\baseX{-2};
\def\xx{1.5};
\def\yy{2};

\node[Path] (i) at (\baseX,\baseHeight) {$i+n$};
\draw[Arr] (i) to +(\xx,0) node[Path, anchor=west](j){$j_1$};

\draw[Arr] (j) to +(\xx,0) node[Path, anchor=west](jn){$j_1+n$};
\draw[Arr] (jn) to +(2*\xx,0) node[Path,anchor=west](k){$k_1$};

\node[Path] (i) at (\baseX, \baseHeight-\yy){$i+n$};
\draw[Arr] (i) to +(1.7*\xx,0) node[Path, anchor=west](j2){$j_2$};
\draw[Arr] (i) to (j2);
\draw[Arr] (j2) to +(\xx,0) node[Path,anchor=west](j2n){$j_2+n$};
\draw[Arr] (j2n) to +(1.3*\xx,0) node[Path,anchor=west]{$k_1$};

\node[Path] (i) at (\baseX,\baseHeight - 2*\yy) {$i+n$};
\draw[Arr] (i) to +(\xx,0) node[Path, anchor=west](j){$j_1$};

\draw[Arr] (j) to +(\xx,0) node[Path, anchor=west](jn){$j_1+n$};
\draw[Arr] (jn) to +(\xx,0) node[Path,anchor=west](k2){$k_2$};

\node[below left = 2.2 cm and 1 cm of i.center,anchor=east](timeStart){};
\def\length{10};
\def\yy{0.4}
\draw[very thick, ->] (timeStart) to +(\length,0) node[below] {$t$};

\path (timeStart) -| node(li){} (i);
\node[below =\yy cm of li.center,anchor=center]{$l_{i+n}$};

\path (timeStart) -| node(ej){} (j); 
\node[below = \yy cm of ej.center,anchor=center]{$e_{j_1}$};

\path (timeStart) -| node(ejn){} (jn); 
\node[below = \yy cm of ejn.center,anchor=center]{$e_{j_1+n}$};

\path (timeStart) -| node(ej2){} (j2); 
\node[below = \yy cm of ej2.center,anchor=center]{$e_{j_2}$};

\path (timeStart) -| node(ej2n){} (j2n); 
\node[below = \yy cm of ej2n.center,anchor=center]{$e_{j_2+n}$};

\path (timeStart) -| node(ek){} (k); 
\node[below = \yy cm of ek.center,anchor=center]{$e_{k_1}$};

\path (timeStart) -| node(ek2){} (k2); 
\node[below = \yy cm of ek2.center,anchor=center]{$e_{k_2}$};

\def\halfLength{0.2};
\foreach \var in {li, ej,ej2,ejn,ej2n,ek,ek2}{
	\draw[thick] (\var.center) to +(0,\halfLength);
	\draw[thick] (\var.center) to +(0,-\halfLength);
}


\end{tikzpicture}

%% file: matchingII.pgf
\begin{tikzpicture}
\usetikzlibrary{arrows,decorations.markings, decorations.pathreplacing, calc}
\usetikzlibrary{math}
\usetikzlibrary{shapes.misc}
\usetikzlibrary{positioning}

\tikzset{Path/.style = {draw,circle, minimum size=1.4cm}}
\tikzset{Arr/.style = {decoration={markings,mark=at position 1 with
    {\arrow[scale=3,>=stealth]{>}}},postaction={decorate}}}
\def\baseHeight{-5.5};
\def\baseX{-2};
\def\xx{4};
\def\yy{4.2};

\node[Path] (i) at (\baseX-0.3, \baseHeight){$i+n$};
\draw[Arr] (i) to node[above]{$2=|\{k_1,k_2\}|$}node[above,yshift=25pt]{\underline{Assignment problem I}} +(\xx,0)node[Path, anchor=west]{$j_1$};
\draw[Arr] (i) to node[above, sloped]{$1=|\{k_1\}|$} +(\xx,-2)node[Path, anchor=west]{$j_2$};

\node[Path, below = \yy cm of i.center,anchor=center] (k2){$k_2$};
\draw[Arr] (k2) to node[above]{$1=|\{i+n\}|$}node[above,yshift=25pt]{\underline{Assignment problem II}} +(\xx,0)node[Path, anchor=west](j1){$j_1$};
\node[below = 2cm of k2.center,anchor=center, Path](k1){$k_1$};
\draw[Arr] (k1) to node[above,sloped]{$1=|\{i+n\}|$} (j1);
\draw[Arr] (k1) to node[above,xshift=3pt]{$1=|\{i+n\}|$} +(\xx,0) node[Path,anchor=west]{$j_2$};

\end{tikzpicture}

%% file: flowchart.pgf
\begin{tikzpicture}[scale=0.95]
\usetikzlibrary{arrows.meta,decorations.markings, decorations.pathreplacing, calc,intersections}
\usetikzlibrary{math}
\usetikzlibrary{shapes.misc}
\usetikzlibrary{positioning}
\tikzset{Box/.style = {draw,rectangle, align=center, minimum width=3.2cm}}
\tikzset{Arr/.style = {-{Latex}}}
\tikzset{  
 connector/.style={
        to path={(\tikztostart) -- ++(#1,0pt) \tikztonodes -| (\tikztotarget) },
		pos=0.5,
		Arr
    },
connector/.default=-2cm}

\def\xx{5};	
\def\yy{0.7};	
\def\yyy{0.6};	

\node[Box] (AA) at (0,0) {LP solution};
\draw[Arr] (AA) to +(0,-\yyy) node[anchor=north,Box](A){Has arc $(i,j)$ been\\ fixed to 1 or 0?};

\path[Arr] (A.south) to +(-\xx,-\yy) node[below, Box] (B) {Merge paths\\ $R_{ki}$ and $R_{jh}$};
\path[Arr] (A.south) to  +(+\xx,-\yy) node[below, Box] (C) {Is there a predecessor\\ relationship between\\ paths $R_{ki}$ and $R_{jh}$?};
\draw[connector] (A.west) to node[above left] {1} (B.north);
\draw[connector=-0cm] (A.east) to (C.north);
\node[above right = 0.01cm and 1cm of A.east] {0};

\draw[Arr] (B.south) to +(0,-\yyy)  node[below, Box] (D) {Is solution\\ infeasible?\\(Subsection \ref{sec:identifyInf})};

\draw[connector=-1cm] (D.west) node[above left]{yes} to +(-1,-3.5) node[below, Box,anchor=north]{Prune node};
\draw[Arr] (D.south) node[below right]{no} to +(0,-\yyy) node [below, Box] (E){Add additional\\arc fixings\\(Subsection \ref{sec:fixVar})};

\draw[Arr](E.east) to node(EF){} +(2,0) node[Box, anchor=west](F){Run phase 1.\\ Have bounds\\been improved?\\(Subsection \ref{sec:generalDesc})};

\draw[connector=-0cm] (C.west)  to (F.north);
\node[above left = 0.01cm and 1.7cm of C.west ] {yes};

\draw[Arr] (F.south) node[below right]{yes} to +(0,-\yyy) node[below, Box] (G) {Run phase 2.\\Have bounds\\been improved?\\(Subsection \ref{sec:generalDesc})};
\draw[connector=-0cm] (G.west) to (EF.center);
\node[above left = 0.9cm and 0.2cm of G.west] {yes};


\path[name path=P1](C.south) to +(0,-5);
\path[name path=P2](F.east) to +(7,0);
\path[name intersections={of=P1 and P2, by=H}];
\node[above=0cm of H.center,anchor=center, Box](H){Add valid\\ inequalities and\\ terminate};
\draw[Arr] (F.east) to node[above]{no} (H.west);
\draw[Arr] (C.south) to node[right]{no} (H.north);
\draw[connector=0cm] (G.east) to (H.south);
\node[above right = 0.01cm and 2cm of G.east] {no};

\end{tikzpicture}

%% file: singlePath.pgf
	
\begin{tikzpicture}
\usetikzlibrary{arrows,decorations.markings, decorations.pathreplacing, calc}
\usetikzlibrary{math}
\usetikzlibrary{shapes.misc}
\usetikzlibrary{positioning}

\tikzset{Path/.style = {draw,circle, minimum size=0.8cm}}
\tikzset{Arr/.style = {decoration={markings,mark=at position 1 with
    {\arrow[scale=3,>=stealth]{>}}},postaction={decorate}}}

\definecolor{myBlue}{RGB}{0,0,0}										
\definecolor{myPurple}{RGB}{0,0,0}									
\definecolor{myMagenta}{RGB}{0,0,0}
\definecolor{myOrange}{RGB}{0,0,0}


\foreach \y in {0}{
	\node[Path, myOrange] (i) at (-4,\y) {$i$};
	\node[Path, myPurple] (j) at (0,\y) {$j$};
	\node[left = 1.5cm of j.center,anchor=center] (jLeft){};	
	\draw[decorate, decoration={brace, amplitude=10pt, raise = 10pt}] (jLeft.west) to node[above=17pt]{Waiting time} (j);

	\draw[Arr] (i) to node[above]{$t_{i,j}$} (jLeft);
}

	\def\lineHeight{-0.8};
	\def\yy{0};
	\def\yUp{0};
	\tikzmath{
		\yy = \lineHeight + 0.15;
		\yUp = (\lineHeight - \yy) * 2;
	}

\def\xx{-4.5};
\draw[->, thick] (\xx,\lineHeight) to (1.5,\lineHeight) node[right] {time}; 
\pgfmathparse{\xx+0.5};
\edef\dummyVal{\pgfmathresult};	
\foreach \x [count=\iii] in {\xx,\dummyVal,...,2}{
}


\node[above =\yy cm of i.center, anchor=center] (ei){}; 
\draw[myOrange](ei.center) to +(0,\yUp) node [below,align=center]{$e_i$\\ $= B_{i}^{LB}$};
\node[right = 1.8cm of ei.center,anchor=center](li){};
\draw[myOrange](li.center) to +(0,\yUp) node[below]{$l_i$};


\node[above left =\yy cm and 0cm of j.center,anchor=center] (ej){};
\draw[myPurple](ej.center) to +(0,\yUp) node [below, align=center]{$e_j$\\$=B_j^{LB}$};
\node[right = 1cm of ej.center,anchor=center](lj){};
\draw[myPurple](lj.center) to +(0,\yUp) node[below]{$l_j$};

\draw[decorate, decoration={brace, amplitude=10pt, raise = 40pt}] (i.center) to node[above=47pt]{Length} (j.center);

\node[align=left] at (-1.5,2.8) {\underline{Path with waiting time}};

\node[] at (-3,-3) (LIJ) {$L_{ij} = B_j^{LB}-B_i^{LB}$};
\node[below = 0.5cm of LIJ.west,anchor=west]{$W_{ij} = B_j^{LB}-(B_i^{LB} + t_{i,j})$};
\node[below = 1cm of LIJ.west,anchor=west]{$T_{ij} = \min(l_i-B_i^{LB},  l_j-B_j^{LB} + W_{ij})$};
\node[above =0.5cm of LIJ.west,anchor=west]{\underline{Path values:}};

\def\yy{-2};
\draw[Arr] (0,\yy-0.5) to +(1,0) node[right]{Fixed arc} ; 

\def\xxx{8};
\draw[Arr] (0+\xxx-0.2,\yy-0.5) to +(1,0) node[right]{Fixed arc} ; 

\pgfmathparse{\xxx-1};
\node[right  = \pgfmathresult cm of i.center,anchor=center](midDummy){};
\draw (midDummy.center) to +(0,-4.3);
\draw (midDummy.center) to +(0,3.3);

\foreach \y in {0}{
	\node[Path, myBlue] (i2) at (-4+\xxx,\y) {$k$};
	\node[Path, myMagenta] (j2) at (0+\xxx,\y) {$h$};

	\draw[Arr] (i2) to node[above]{$t_{k,h}$} (j2);
}

	\def\lineHeight{-0.8};
	\def\yy{0};
	\def\yUp{0};
	\tikzmath{
		\yy = \lineHeight + 0.15;
		\yUp = (\lineHeight - \yy) * 2;
	}

\def\xx{-4.5};
\draw[->, thick] (\xx+\xxx,\lineHeight) to (1.5+\xxx,\lineHeight) node[right] {time}; 
\pgfmathparse{\xx+0.5};
\edef\dummyVal{\pgfmathresult};	
\foreach \x [count=\iii] in {\xx,\dummyVal,...,2}{
}


\node[above =\yy cm of i2.center, anchor=center] (ei){}; 
\draw[myBlue](ei.center) to +(0,\yUp) node [below,align=center]{$e_k$\\ $= B_{k}^{LB}$};
\node[right = 1.8cm of ei.center,anchor=center](li){};
\draw[myBlue](li.center) to +(0,\yUp) node[below]{$l_k$};


\node[above left =\yy cm and 0cm of j2.center,anchor=center] (Bj){};
\draw[myMagenta](Bj.center) to +(0,\yUp) node [below, align=center]{$B_h^{LB}$};
\node[left = 1cm of Bj.center, anchor=center](ej){};
\draw[myMagenta](ej.center) to +(0,\yUp) node [below]{$e_h$};
\node[right = 0.8cm of Bj.center,anchor=center](lj){};
\draw[myMagenta](lj.center) to +(0,\yUp) node[below]{$l_h$};

\draw[decorate, decoration={brace, amplitude=10pt, raise = 25pt}] (i2.center) to node[above=32pt]{Length} (j2.center);

\node[align=left] at (-1.5+\xxx,2.8) {\underline{Path without waiting time}};

\node[] at (-3+\xxx,-3) (LIJ) {$L_{kh} = B_h^{LB}-B_k^{LB}$};
\node[below = 0.5cm of LIJ.west,anchor=west]{$W_{kh} = 0$};
\node[below = 1cm of LIJ.west,anchor=west]{$T_{kh} = \min(l_k-B_k^{LB},  l_h-B_h^{LB})$};
\node[above =0.5cm of LIJ.west,anchor=west]{\underline{Path values:}};

\end{tikzpicture}

%% file: phaseOne.pgf
\begin{tikzpicture}
\usetikzlibrary{arrows,decorations.markings, calc, decorations.pathreplacing}
\usetikzlibrary{math}
\usetikzlibrary{shapes.misc}
\usetikzlibrary{positioning}

	\tikzset{Path/.style = {draw,circle}}
	\tikzset{Arr/.style = {decoration={markings,mark=at position 1 with
    	{\arrow[scale=3,>=stealth]{>}}},postaction={decorate}}}

	\node[Path] (P1) at (0,0) {$P_1$};
	\node[Path, right = 4cm of P1.center,anchor=center] (P2) {$P_2$};
	\draw[Arr] (P1) to (P2);
	\draw[decorate, decoration={brace, amplitude=10pt, raise = 15pt}] (P1.west) to node[above=22pt]{Merged path} (P2.east);
	
	\node[Path, above right = 1.5cm and 2cm of P2.center, anchor=center, label = {\scriptsize $B^{LB}_{D_1}\geq B^{LB}_{P_2} + tt_{P_2,D_1}$}] (D1) {$D_1$};	
	\node[Path, right = 3cm of D1.center, anchor=center] (P3){$P_3$};	
	\draw[Arr, dashed] (P2) to (D1);
	\draw[Arr] (D1) to (P3);
	
	\node[Path, below right = 1.5cm and 2cm of P2.center,anchor=center,label = 270:{\scriptsize $B^{LB}_{D_2}\geq B^{LB}_{P_2} + tt_{P_2,D_2}$}] (D2) {$D_2$};
	\draw[Arr, dashed] (P2) to (D2); 

	\node[above right = 2cm and 2cm of P3.center, anchor=center] (D3) {$\cdots$};
	\draw[Arr, dashed] (P3) to (D3);

	\def\lineHeight{-2.8};
	\def \yy{0};
	\def\yUp{0};
	\tikzmath{
		\yy = \lineHeight + 0.15;
		\yUp = (\lineHeight - \yy) * 2;
	}
	
	\draw[->, thick] (-1,\lineHeight) to (10,\lineHeight) node[right] {time}; 

	\node[above = \yy cm of P1.center,anchor=center] (BP1) {};
	\draw (BP1.center) to +(0,\yUp) node[below]{$B_{P_1}^{LB}$};

	\node[above = \yy cm of P2.center,anchor=center] (BP2){};
	\draw (BP2.center) to +(0,\yUp) node[below,align=center]{$B_{P_2}^{LB} =$\\$B_{P_1}^{LB} + L_{P_1,P_2}=$\\ $\max(e_{P_2}, B_{P_1}^{LB} + t_{P_1,P_2})$};

	\def\xx{5.5};
	\def\yy{0};
	\draw[Arr, dashed] (\xx,\yy) to +(1,0) node[right]{Shortest possible connection} ; 
	\draw[Arr] (\xx,\yy-0.5) to +(1,0) node[right]{Fixed arc} ; 
\end{tikzpicture}

%% file: pathResequence.pgf
\begin{tikzpicture}
\usetikzlibrary{arrows.meta,decorations.markings, calc, decorations.pathreplacing}
\usetikzlibrary{math}
\usetikzlibrary{shapes.misc}
\usetikzlibrary{positioning}

	\tikzset{Path/.style = {draw,circle}}
	\tikzset{Arr/.style = {decoration={markings,mark=at position 1 with
    	{\arrow[scale=3,>=stealth]{>}}},postaction={decorate}}}

\def\yOne{2}	
\def\xOne{2.5};	
\def\yTwo{0.5}; 
\foreach \y in {1,2,3}{
	\node[Path] (P\y) at (0,-\y*\yOne) {$P_\y$};

}
\node[above right = \yTwo cm and \xOne cm of P1.center, anchor=center,Path ](P12){$P_2$};
\node[below right = \yTwo cm and \xOne cm of P1.center, anchor=center,Path ](P13){$P_3$};
\node[above right = \yTwo cm and \xOne cm of P2.center, anchor=center,Path ](P21){$P_1$};
\node[below right = \yTwo cm and \xOne cm of P2.center, anchor=center,Path ](P23){$P_3$};
\node[above right = \yTwo cm and \xOne cm of P3.center, anchor=center,Path ](P31){$P_1$};
\node[below right = \yTwo cm and \xOne cm of P3.center, anchor=center,Path ](P32){$P_2$};

\node[right = 2*\xOne cm of P1.center,anchor=center,Path](P33) {$P_3$};
\node[right = 2*\xOne cm of P2.center,anchor=center,Path](P11) {$P_1$};
\node[right = 2*\xOne cm of P3.center,anchor=center,Path](P22) {$P_2$};

\node[right = 3*\xOne cm of P2.center,anchor=center,Path](D1){$D_{}$};	

\node[above right = 0.7cm and 1.2cm of P33.center,anchor=center]{\scriptsize$S=\{P_1, P_2\}, v=\{P_3\}$};

\node[below left = 0.6*\yOne cm and 0.5cm of P3.center,anchor=center](St){};
\draw[-Latex](St) to +(3.5*\xOne,0)node[below]{Stages};
\node[below = 0.7cm of St.center,anchor=center, Path](lg){$P$};
\node[right = 1.6cm of lg.center,anchor=center]{\scriptsize State with $v=P$};

\foreach \y in {1,2,3}{	
	\foreach \x in {1,2,3}{
		\ifnum \y=\x 
		[\else
			\draw(P\y) to (P\y\x);
			\foreach \xx in {1,2,3}{
				\ifnum \x=\xx
				[\else
						\ifnum\y=\xx
							[\else
								\draw(P\y\x) to (P\xx\xx);
							]\fi
				]\fi
			}

		 ]\fi
	}
	\draw(P\y\y) to (D1);
}

\end{tikzpicture}

%% file: DARP_MIP_TW.bbl
\begin{thebibliography}{}

\bibitem[\protect\astroncite{Baldacci et~al.}{2011}]{BBM11}
Baldacci, R., Bartolini, E., and Mingozzi, A. (2011).
\newblock An exact algorithm for the pickup and delivery problem with time
  windows.
\newblock {\em Operations Research}, 59(2):414--426.

\bibitem[\protect\astroncite{Bellman}{1962}]{Bel62}
Bellman, R. (1962).
\newblock Dynamic programming treatment of the travelling salesman problem.
\newblock {\em Journal of the ACM (JACM)}, 9(1):61--63.

\bibitem[\protect\astroncite{Braekers et~al.}{2014}]{BCJ14}
Braekers, K., Caris, A., and Janssens, G.~K. (2014).
\newblock Exact and meta-heuristic approach for a general heterogeneous
  dial-a-ride problem with multiple depots.
\newblock {\em Transportation Research Part B: Methodological}, 67:166--186.

\bibitem[\protect\astroncite{Camm et~al.}{1990}]{CRT90}
Camm, J.~D., Raturi, A.~S., and Tsubakitani, S. (1990).
\newblock Cutting big {M} down to size.
\newblock {\em Interfaces}, 20(5):61--66.

\bibitem[\protect\astroncite{Chevrier et~al.}{2012}]{CLJD12}
Chevrier, R., Liefooghe, A., Jourdan, L., and Dhaenens, C. (2012).
\newblock Solving a dial-a-ride problem with a hybrid multi-objective
  evolutionary approach: Application to demand responsive transport.
\newblock {\em Applied Soft Computing}, 12(4):1247--1258.

\bibitem[\protect\astroncite{Codato and Fischetti}{2006}]{CF06}
Codato, G. and Fischetti, M. (2006).
\newblock Combinatorial {B}enders' cuts for mixed-integer linear programming.
\newblock {\em Operations Research}, 54(4):756--766.

\bibitem[\protect\astroncite{Cordeau}{2006}]{Cor06}
Cordeau, J.-F. (2006).
\newblock A branch-and-cut algorithm for the dial-a-ride problem.
\newblock {\em Operations Research}, 54(3):573--586.

\bibitem[\protect\astroncite{Cort{\'e}s et~al.}{2010}]{CMC10}
Cort{\'e}s, C.~E., Matamala, M., and Contardo, C. (2010).
\newblock The pickup and delivery problem with transfers: Formulation and a
  branch-and-cut solution method.
\newblock {\em European Journal of Operational Research}, 200(3):711--724.

\bibitem[\protect\astroncite{Diana and Dessouky}{2004}]{DD04}
Diana, M. and Dessouky, M.~M. (2004).
\newblock A new regret insertion heuristic for solving large-scale dial-a-ride
  problems with time windows.
\newblock {\em Transportation Research Part B: Methodological}, 38(6):539--557.

\bibitem[\protect\astroncite{Dijkstra}{1959}]{D59}
Dijkstra, E.~W. (1959).
\newblock A note on two problems in connexion with graphs.
\newblock {\em Numerische Mathematik}, 1(1):269--271.

\bibitem[\protect\astroncite{Dror}{1994}]{D94}
Dror, M. (1994).
\newblock Note on the complexity of the shortest path models for column
  generation in {VRPTW}.
\newblock {\em Operations Research}, 42(5):977--978.

\bibitem[\protect\astroncite{Gschwind and Drexl}{2019}]{GD19}
Gschwind, T. and Drexl, M. (2019).
\newblock Adaptive large neighborhood search with a constant-time feasibility
  test for the dial-a-ride problem.
\newblock {\em Transportation Science}, 53(2):480--491.

\bibitem[\protect\astroncite{Gschwind and Irnich}{2015}]{GI15}
Gschwind, T. and Irnich, S. (2015).
\newblock Effective handling of dynamic time windows and its applications to
  solving the dial-a-ride problem.
\newblock {\em Transportation Science}, 49(2):335--354.

\bibitem[\protect\astroncite{Held and Karp}{1962}]{HK62}
Held, M. and Karp, R.~M. (1962).
\newblock A dynamic programming approach to sequencing problems.
\newblock {\em Journal of the Society for Industrial and Applied Mathematics},
  10(1):196--210.

\bibitem[\protect\astroncite{Ho et~al.}{2018}]{HSKLPT18}
Ho, S.~C., Szeto, W.~Y., Kuo, Y.-H., Leung, J. M.~Y., Petering, M., and Tou, T.
  W.~H. (2018).
\newblock A survey of dial-a-ride problems: Literature review and recent
  developments.
\newblock {\em Transportation Research Part B: Methodological}, 111:395--421.

\bibitem[\protect\astroncite{Ioachim et~al.}{1998}]{IGSD98}
Ioachim, I., Glinas, S., Soumis, F., and Desrosiers, J. (1998).
\newblock A dynamic programming algorithm for the shortest path problem with
  time windows and linear node costs.
\newblock {\em Networks}, 31(3):193--204.

\bibitem[\protect\astroncite{Kalantari et~al.}{1985}]{KHA85}
Kalantari, B., Hill, A.~V., and Arora, S.~R. (1985).
\newblock An algorithm for the traveling salesman problem with pickup and
  delivery customers.
\newblock {\em European Journal of Operational Research}, 22(3):377--386.

\bibitem[\protect\astroncite{Little et~al.}{1963}]{LMSK63}
Little, J. D.~C., Murty, K.~G., Sweeney, D.~W., and Karel, C. (1963).
\newblock An algorithm for the traveling salesman problem.
\newblock {\em Operations Research}, 11(6):972--989.

\bibitem[\protect\astroncite{Lu and Dessouky}{2004}]{LD04}
Lu, Q. and Dessouky, M.~M. (2004).
\newblock An exact algorithm for the multiple vehicle pickup and delivery
  problem.
\newblock {\em Transportation Science}, 38(4):503--514.

\bibitem[\protect\astroncite{Lysgaard}{2006}]{L06}
Lysgaard, J. (2006).
\newblock Reachability cuts for the vehicle routing problem with time windows.
\newblock {\em European Journal of Operational Research}, 175(1):210--223.

\bibitem[\protect\astroncite{Molenbruch et~al.}{2017}]{MBCV17}
Molenbruch, Y., Braekers, K., Caris, A., and {Vanden Berghe}, G. (2017).
\newblock Multi-directional local search for a bi-objective dial-a-ride problem
  in patient transportation.
\newblock {\em Computers {\&} Operations Research}, 77:58--71.

\bibitem[\protect\astroncite{Paquette et~al.}{2013}]{PCLP13}
Paquette, J., Cordeau, J.-F., Laporte, G., and Pascoal, M.~M. (2013).
\newblock Combining multicriteria analysis and tabu search for dial-a-ride
  problems.
\newblock {\em Transportation Research Part B: Methodological}, 52:1--16.

\bibitem[\protect\astroncite{Parragh et~al.}{2008}]{PDH08}
Parragh, S.~N., Doerner, K.~F., and Hartl, R.~F. (2008).
\newblock A survey on pickup and delivery models: Part ii: Transportation
  between pickup and delivery locations.
\newblock {\em Journal f{\"u}r Betriebswirtschaft}, 58(2):81--117.

\bibitem[\protect\astroncite{Parragh et~al.}{2009}]{PDHG09}
Parragh, S.~N., Doerner, K.~F., Hartl, R.~F., and Gandibleux, X. (2009).
\newblock A heuristic two-phase solution approach for the multi-objective
  dial-a-ride problem.
\newblock {\em Networks}, 54(4):227--242.

\bibitem[\protect\astroncite{Parragh et~al.}{2015}]{PPA15}
Parragh, S.~N., {Pinho de Sousa}, J., and Almada-Lobo, B. (2015).
\newblock The dial-a-ride problem with split requests and profits.
\newblock {\em Transportation Science}, 49(2):311--334.

\bibitem[\protect\astroncite{Pfeiffer and Schulz}{2022a}]{PS20}
Pfeiffer, C. and Schulz, A. (2022a).
\newblock An {ALNS} algorithm for the static dial-a-ride problem with ride and
  waiting time minimization.
\newblock {\em OR Spectrum}, 44:87--119.

\bibitem[\protect\astroncite{Pfeiffer and Schulz}{2022b}]{PS22}
Pfeiffer, C. and Schulz, A. (2022b).
\newblock A new lower bound for the static dial-a-ride problem with ride and
  waiting time minimization.
\newblock In Freitag, M., Kinra, A., Kotzab, H., and Megow, N., editors, {\em
  Dynamics in Logistics}, Lecture Notes in Logistics, pages 231--243. {Springer
  International Publishing}.

\bibitem[\protect\astroncite{Psaraftis}{1980}]{Psa80}
Psaraftis, H.~N. (1980).
\newblock A dynamic programming solution to the single vehicle many-to-many
  immediate request dial-a-ride problem.
\newblock {\em Transportation Science}, 14(2):130--154.

\bibitem[\protect\astroncite{Psaraftis}{1983}]{Psa83}
Psaraftis, H.~N. (1983).
\newblock An exact algorithm for the single vehicle many-to-many dial-a-ride
  problem with time windows.
\newblock {\em Transportation Science}, 17(3):351--357.

\bibitem[\protect\astroncite{Qu and Bard}{2015}]{QB15}
Qu, Y. and Bard, J.~F. (2015).
\newblock A branch-and-price-and-cut algorithm for heterogeneous pickup and
  delivery problems with configurable vehicle capacity.
\newblock {\em Transportation Science}, 49(2):254--270.

\bibitem[\protect\astroncite{Riedler and Raidl}{2018}]{RR18}
Riedler, M. and Raidl, G. (2018).
\newblock Solving a selective dial-a-ride problem with logic-based benders
  decomposition.
\newblock {\em Computers {\&} Operations Research}, 96:30--54.

\bibitem[\protect\astroncite{Rist and Forbes}{2022}]{RF22}
Rist, Y. and Forbes, M. (2022).
\newblock A column generation and combinatorial benders decomposition algorithm
  for the selective dial-a-ride-problem.
\newblock {\em Computers {\&} Operations Research}, 140:105649.

\bibitem[\protect\astroncite{Rist and Forbes}{2021}]{RF21}
Rist, Y. and Forbes, M.~A. (2021).
\newblock A new formulation for the dial-a-ride problem.
\newblock {\em Transportation Science}, 55(5):1113--1135.

\bibitem[\protect\astroncite{Ropke and Cordeau}{2009}]{RC09}
Ropke, S. and Cordeau, J.-F. (2009).
\newblock Branch and cut and price for the pickup and delivery problem with
  time windows.
\newblock {\em Transportation Science}, 43(3):267--286.

\bibitem[\protect\astroncite{Ropke et~al.}{2007}]{RCL07}
Ropke, S., Cordeau, J.-F., and Laporte, G. (2007).
\newblock Models and branch-and-cut algorithms for pickup and delivery problems
  with time windows.
\newblock {\em Networks}, 49(4):258--272.

\bibitem[\protect\astroncite{Ropke and Pisinger}{2006}]{RP06}
Ropke, S. and Pisinger, D. (2006).
\newblock An adaptive large neighborhood search heuristic for the pickup and
  delivery problem with time windows.
\newblock {\em Transportation Science}, 40(4):455--472.

\bibitem[\protect\astroncite{Ruland and Rodin}{1997}]{RR97}
Ruland, K.~S. and Rodin, E.~Y. (1997).
\newblock The pickup and delivery problem: Faces and branch-and-cut algorithm.
\newblock {\em Computers {\&} Mathematics with Applications}, 33(12):1--13.

\bibitem[\protect\astroncite{Savelsbergh}{1992}]{S92}
Savelsbergh, M.~W. (1992).
\newblock The vehicle routing problem with time windows: Minimizing route
  duration.
\newblock {\em ORSA Journal on Computing}, 4(2):146--154.

\bibitem[\protect\astroncite{Tachet et~al.}{2017}]{TSSRSSR17}
Tachet, R., Sagarra, O., Santi, P., Resta, G., Szell, M., Strogatz, S.~H., and
  Ratti, C. (2017).
\newblock Scaling law of urban ride sharing.
\newblock {\em Scientific Reports}, 7(1):1--6.

\end{thebibliography}
